\documentclass[reqno]{article}
\usepackage{amssymb, amsmath, amsthm, amsfonts, amscd}
\usepackage[english]{babel}
\usepackage{bbm}

\usepackage{graphicx}
\graphicspath{ {./images/} }
\usepackage{epstopdf}
\usepackage[margin=3.cm]{geometry}

\usepackage{breqn}
\usepackage{mathtools}
\usepackage{color}
\usepackage{hyperref}
\hypersetup{
    colorlinks=true,
    citecolor=red,
    filecolor=black,
    linkcolor=blue,
    urlcolor=blue
}

\usepackage{siunitx}
\usepackage{booktabs}

\usepackage{xr}
\chardef\bslash=`\\ 

\newtheorem{theorem}{Theorem}[section]
\newtheorem*{theorem*}{Theorem}
\newtheorem{corollary}[theorem]{Corollary}
\newtheorem{lemma}[theorem]{Lemma}
\newtheorem{prop}[theorem]{Proposition}

\newtheorem*{remark}{Remark}

\newcommand{\N}{\mathbb{N}}

\newcommand{\B}{\mathrm{B}}
\newcommand{\F}{\mathrm{F}}
\newcommand{\FF}{\mathbb{F}}
\newcommand{\FFF}{\mathcal{F}}

\newcommand{\PP}{\mathcal{P}}

\newcommand{\R}{\mathbb{R}}

\newcommand{\T}{\bar{T}}

\newcommand{\HH}{H^{1}}
\newcommand{\dx}{\dot{x}}
\newcommand{\h}{^{0}_{0}}
\newcommand{\dd}{\mathrm{d}}

\def\m{\mu }
\def\n{\nu }

\def\s{\sigma}
\def\t{\tau}

\def\d{\delta}

\def\w{\omega}

\def\f{\varphi}
\def\.{\cdot }
\def\ra{\rightarrow}

\title{
\textsc{\textbf{The weak form of the SDOF and MDOF equation of motion, part II: A numerical method for the SDOF problem}}\\
\author{Nikolaos Karaliolios$^1$ \and Dimitrios L. Karabalis$^2$}
\date{%
    $^1$email: nkaraliolios@gmail.com\\%
    $^2$Department of Civil Engineering, University of Patras\\[2ex]%
    \today
}
}

\begin{document}

\maketitle

\begin{abstract}
A new, more efficient, numerical method for the SDOF problem is presented.
Its construction is based on the weak form of the equation of motion, as obtained
in \cite{NKSDOFI}, using piece-wise polynomial functions as interpolation
functions. The approximation rate can be arbitrarily high, proportional to the degree
of the interpolation functions, tempered only by numerical instability. Moreover,
the mechanical energy of the system is conserved. Consequently, all significant
drawbacks of existing algorithms, such as the limitations imposed by the Dahlqvist
Barrier theorem \cite{Dahl56} and the need for introduction of numerical damping, have
been overcome.
\end{abstract}  

\tableofcontents

\section{Application of the weak formulation to piece-wise polynomial interpolation} \label{secnum method}

The algorithm presented herein, based on the weak formulation of the SDOF problem as
obtained in \cite{NKSDOFI}, is constructed as follows. A partition of
the interval $[0,\T ]$ into $l$ intervals of length $h$ is considered, and the approximate solution $x _{ap}$ will be polynomial in each
subinterval of the type
\begin{equation}
I_{j} = [jh,(j+1)h]
\end{equation}
The space formed by such functions will be denoted by
$\PP _{h,p} = \PP _{h,p}([0,\T ])$.

The breakpoints $\{ j h \}_{j = 0 \cdots l}$, i.e. the points where the expression of a
function in $\PP _{h,p}$ is allowed to change are equally spaced in $I$, and the
corresponding time-step is $h > 0$. The time-step $h$ is chosen so that
\begin{equation}
\frac{\T}{h} = l \in \N ^{*}
\end{equation}
where $l$ is the number of time-steps of the algorithm.

The functions $B_{i,p,h}$, defined here below, for some fixed $p\geq 3$, will form a basis
for the restriction of functions in $\PP _{h,p}$ in each subinterval formed by consecutive
break-points. The weak formulation of the SDOF problem will be considered in each such
interval $I_{j}$. The algorithm thus takes automatically a time-step
algorithm character and, even though it does not represent the full potential of the weak
formulation, can be seen to be already more powerful than those in the literature.

Consequently, there is some change in notation with respect to \cite{NKSDOFI}, since $\T$ still
stands for the upper bound of the interval in which the solution is approximated, but
the intervals on which the weak formulation is considered are precisely the $I_{j}$, for
$j = 0, \cdots , l-1$ where $l = \T / h \in \N ^{*}$.

The authors are already working on the construction of a numerical method based on the
weak formulation on the entire interval $[0,\T ]$, which will also have a step-by-step
nature.

\section{Choice of functions}

Fix $h >0$. The space
\begin{equation}
\PP= \PP(p,h) \subset \HH ([0,h])
\end{equation}
is chosen to be equal to the space of polynomial functions of a given degree
$p -1\in \N ^{*}$, $p \geq 3$. It is of dimension $p$.
The space
\begin{equation}
\dot{\PP}= \dot{\PP} (p,h)= \{ \dot{\f}, \f \in
\PP(p,h) 
\}
\subset \HH ([0,h])
\end{equation}
is thus equal to the space of polynomial functions of degree $p -2\in \N ^{*}$ and of
dimension $p-1$.

The basis of $\PP $ given in the next equation is known in the literature as the
Bernstein polynomials, see, e.g., \cite{Kac1938}. Call, for $i=1,\cdots , p$,
\begin{equation}
B_{i,p,h}(t) = B_{i}(t) = \begin{cases}
\frac{\Gamma (p)}{\Gamma (i) \Gamma (p-i +1)}
\left( \frac{t}{h}
\right)^{i-1}
\left( \frac{h-t}{h}
\right)^{p-i} \text{ if } 0 \leq t \leq h \\
0 \hphantom{
\frac{\Gamma (p)}{\Gamma (i) -}
\left( \frac{t}{h}
\right)^{i-1}
} \text{   otherwise}
\end{cases}
\end{equation}
All functions $B_{i}$ are by their definition non-negative.
The normalization factor
\begin{equation}
\frac{\Gamma (p)}{\Gamma (i) \Gamma (p-i +1)} = {{p-1}\choose {i-1}}
=
\frac{(p-1)!}{(i-1)!(p-i)!}
\end{equation}
is chosen so that
\begin{equation}
\sum _{i=1}^{p}
B_{i}(t) = \begin{cases}
1 \text{ if } 0 \leq t \leq h \\
0  \text{   otherwise}
\end{cases}
\end{equation}
i.e. so that the $B_{i}$ form a \textit{partition of unity} on the interval $[0,h]$ (cf. \cite{dBSpl}, \S IX).





\section{Some useful properties of Bernstein polynomials} \label{secprop bernstein}

\subsection{The functions and their derivatives}

The Bernstein polynomials form a basis of polynomial functions of degree $p-1$ defined on
$[0,h]$. For this particular choice of basis, one can immediately verify that
\begin{equation}
\begin{array}{l@{}l}
B_{i}(0) \neq 0& \iff i = 1 \\
B_{i}(h) \neq 0 &\iff i = p \\
B_{i}'(0) \neq 0 & \iff i \in \{ 1,2 \} \\
B_{i}'(h) \neq 0 &\iff i \in \{ p-1,p \} \\
\end{array}
\end{equation}
More precisely,
\begin{equation}
\begin{array}{l@{}l}
B_{1}(0) = B_{p}(h) &=1 \\
B_{1}'(0) = - B_{p}'(h) &=-\frac{p-1}{h} \\
B_{2}(0) = B_{p-1}(h) &=0 \\
B_{2}'(0) = - B_{p-1}'(h) &=\frac{p-1}{h}
\end{array}
\end{equation}

The following lemma will be useful for the calculation of the derivatives of Bernstein
polynomials. The derivative of a polynomial of degree $p$ is a polynomial of degree
$p-1$. The following lemma shows how to express a Bernstein polynomial of degree $p-1$
in terms of Bernstein polynomials of degree $p$.
\begin{lemma} \label{lemlower to higher}
The following formula holds true for $i = 1,\cdots , p-1$:
\begin{equation}
B_{i,p-1,h}(t) =
\frac{p-i}{p-1}
B_{i,p,h}(t) + 
\frac{i}{p-1}
B_{i+1,p,h}(t)
\end{equation}
\end{lemma}
\begin{proof}
The proof is elementary and is based on the following observation:
\begin{equation}
\begin{array}{r@{}l}
\left( \frac{t}{h}
\right)^{i-1}
\left( \frac{h-t}{h}
\right)^{p-i}
+
\left( \frac{t}{h}
\right)^{i}
\left( \frac{h-t}{h}
\right)^{p-i-1}
&=
\left( \frac{t}{h}
\right)^{i-1}
\left( \frac{h-t}{h}
\right)^{p-i-1}
\left(
\frac{t}{h}
+
\frac{h-t}{h}
\right)
\\
&=
\left( \frac{t}{h}
\right)^{i-1}
\left( \frac{h-t}{h}
\right)^{p-i-1}
\\
&=
{{p-2}\choose {i-1}}^{-1}B_{i,p-1,h}(t)
\end{array}
\end{equation}
and the rest is a direct calculation involving the normalization factors.
\end{proof}

The following lemma is elementary. The proof is by immediate calculation and
use of lemma \ref{lemlower to higher}.
\begin{lemma} \label{lemBern deriv}
The following formula holds true for $i = 1,\cdots , p$, where for convenience,
for all $p \in \N ^{*} $, $B_{0,p,h}(t) \equiv B_{p+1,p,h}(t) \equiv 0$:
\begin{equation}
\begin{array}{r@{}l}
\frac{d}{d t}
B_{i,p,h}(t) &=
\frac{p-1}{h}
(
B_{i-1	,p-1,h}(t) -
B_{i,p-1,h}(t)
)
\\
&=
\frac{1}{h}
(
(p-i-1)
B_{i-1,p,h}(t) + 
(2i-p-1)
B_{i,p,h}(t)-
i
B_{i+1,p,h}(t)
)
\end{array}
\end{equation}
\end{lemma}

The following corollary follows from direct calculation.
\begin{corollary}
The local maximum of $B_{i,p,h}$ for $i=2,\cdots,p-1$ is attained at $t= \frac{ih}{p}$
\end{corollary}

\begin{proof}
From the formula of the lemma it easily follows that
\begin{equation}
\frac{d}{d t}
B_{i,p,h}(t) =
\frac{(p-1)!}{i!(p-i)!}
\left(
\frac{t}{h}
\right)^{i-1}
\left(
\frac{h-t}{h}
\right)^{p-i-1}
\frac{ih-pt}{h}
\end{equation}
whose only root inside $(0,h)$ is $t=\frac{ih}{p}$ and the corollary follows.
\end{proof}

\subsubsection{Integrals and $L^{2}$ products}

The following lemma concerns the integral of a Bernstein polynomial. It will be useful
for the calculation of $L^{2}_{c}$ scalar products.
\begin{lemma} \label{lemBern integral}
The integral of a Bernstein polynomial $B_{i,p,h}$ over $[0,h]$ depends only on $p$ and
$h$:
\begin{equation}
\int _{0}^{h}
B_{i,p,h}(t) \dd t
=\frac{h}{p}
\end{equation}
\end{lemma}

\begin{proof}
The proof can proceed by integration by parts. Alternatively, lemma \ref{lemBern deriv}
implies that, for $i = 1,\cdots , p$,
\begin{equation}
\begin{array}{r@{}l}
\int _{0}^{h}
(
B_{i-1	,p,h}(t) -
B_{i,p,h}(t)
)dt
&=
\frac{h}{p}
\int _{0}^{h}
\frac{d}{d t}
B_{i,p+1,h}(t) \dd t
\\
&=
\frac{h}{p}
(
B_{i,p+1,h}(h)
-
B_{i,p+1,h}(0))
\\
&=
\frac{h}{p}
(
\d _{i,p+1}
-
\d _{i,1}
)
\end{array}
\end{equation}
where $\d _{kl}$ is the Kronecker delta function, equal to $1$ if $k=l$ and $0$ otherwise.
The conclusion of the lemma follows directly by using the fact that the Bernstein
polynomials form a partition of unity, which implies that
\begin{equation}
\sum _{i=1}^{p}
\int _{0}^{h}
B_{i,p,h}(t) \dd t
= h
\end{equation}
from which the conclusion follows easily.
\end{proof}

A formula for indefinite integrals of Bernstein polynomials can also be obtained.
The proof is by integration by parts and application of lemma \ref{lemBer momentum}
for $q = 1$.
\begin{lemma} \label{lemBern integrals}
It holds true that
\begin{equation}
\begin{array}{r@{}l}
    \int _{0}^{t} B_{i,p,h}(\t )\mathrm{d} \t &=
    \frac{h}{p}
    \sum _{j=i+1}^{p+1}
    B_{j,p+1,h}(t )
    \\
    &=
    \frac{h}{p}\left(1-
    \sum _{j=1}^{i}
    B_{j,p+1,h}(t )
    \right)
\end{array}
\end{equation}
\end{lemma}

\begin{proof}
    Direct calculation gives
    \begin{equation}
\begin{array}{r@{}l}
\int _{0}^{t}  B_{i,p,h}(\tau )\mathrm{d}\tau
&=
tB_{i,p,h}(t ) -
\int _{0}^{t} \tau \dot{B}_{i,p,h}(\tau )\mathrm{d}\tau
\\
&=
h\frac{i}{p}
B_{i+1,p+1,h}(t ) -
\frac{p-1}{h}
\int _{0}^{t} \tau (B_{i-1,p-1,h}(\tau ) 
- B_{i,p-1,h}(\tau ))
\mathrm{d}\tau
\\
&=
h\frac{i}{p}
B_{i+1,p+1,h}(t ) -
(i-1)
\int _{0}^{t} B_{i,p,h}(\tau ) 
\mathrm{d}\tau
+
i
\int _{0}^{t} \tau B_{i+1,p,h}(\tau ) 
\mathrm{d}\tau
\end{array}
\end{equation}
so that, all in all,
\begin{equation}
    \int _{0}^{t}  B_{i,p,h}(\tau )\mathrm{d}\tau
    =
    \frac{h}{p}
    B_{i+1,p+1,h}(t ) +
    \int _{0}^{t}  B_{i+1,p,h}(\tau )\mathrm{d}\tau
\end{equation}
It can be seen directly, either by calculation or by the
standing convention that $B_{p+1,p,h} \equiv 0$, that
\begin{equation}
    \int _{0}^{t}  B_{p,p,h}(\tau )\mathrm{d}\tau
    =
    \frac{h}{p}
    B_{p+1,p+1,h}(t )
\end{equation}
so that, for all $1\leq i \leq p$,
\begin{equation}
    \int _{0}^{t}  B_{i,p,h}(\tau )\mathrm{d}\tau
    =
    \frac{h}{p}
    \sum _{j=i+1}^{p+1}
    B_{j,p+1,h}(t )
\end{equation}
In particular, lemma \ref{lemBern integral} can be obtained
by this formula.

The last part follows from the fact that the Bernstein polynomials
form a partition of unity.
\end{proof}

The product of two Bernstein polynomials satisfies the following simple rule, the proof
of which is left to the reader.
\begin{lemma} \label{lemBern prod}
It holds true that, for $1\leq i \leq p$ and $1 \leq j \leq q$ and $h>0$,
\begin{equation}
B_{i,p,h}(t) B_{j,q,h}(t) = \frac{{{p-1}\choose{i-1}}{{q-1}\choose{j-1}}}{{{p+q-2}\choose{i+j-2}}}
B_{i+j-1,p+q-1,h}(t)
\end{equation}
\end{lemma}

Combining lemmas \ref{lemBern integral} and \ref{lemBern prod} directly yields the
following corollary. For notation, the reader is referred to eq.
(22) of \cite{NKSDOFI}.
\begin{corollary} \label{corBern prod L2}
The $L^{2}_{0}$ scalar product of two Bernstein polynomials is given by
\begin{equation}
\langle B_{i,p,h}(\. ), B_{j,q,h}(\. ) \rangle_{0}=
\int _{0}^{h} B_{i,p,h}(\. ) B_{j,q,h}(\. ) =
\frac{{{p-1}\choose{i-1}}{{q-1}\choose{j-1}}}{(p+q-1){{p+q-2}\choose{i+j-2}}}
h
\end{equation}
The $L^{2}_{0}$ norm of a Bernstein polynomial is given by
\begin{equation}
\| B_{i,p,h}(\. ) \|_{L^{2}_{0}} =
\left(
\langle B_{i,p,h}(\. ), B_{i,p,h}(\. ) \rangle_{0}
\right)^{1/2} =
\frac{{{p-1}\choose{i-1}}}{\sqrt{(2p-1){{2p-2}\choose{2i-2}}}}
\sqrt{h}
\end{equation}
\end{corollary}

A numerical estimation of the behavior of products is given in the following
lemma.
\begin{lemma} \label{lemBern prod numerical}
For $3 \leq p \leq 200$ it holds true that
\begin{equation}
\begin{array}{r@{}l}
\max \limits _{1\leq i,j\leq p}
\langle B_{i,p,h} (\.), B_{j,p,h} (\.) \rangle_{0} &\leq
\frac{h}{p^{11/10}}
\end{array}
\end{equation}
\end{lemma}

In particular, since
\begin{equation}
t^{q} = h^{q}B_{q+1,q+1,h}(t) \text{ for } t \in [0,h]
\end{equation}
one gets directly the following lemma, concerning the moments of the Bernstein polynomials.
It will be useful in the estimation of $L^{2}_{c}$ scalar products, using the formula
for the $L^{2}_{0}$ scalar product.
\begin{lemma} \label{lemBer momentum}
For $1\leq i \leq p$  and $h>0$,
\begin{equation}
\begin{array}{r@{}l}
t^{q}B_{i,p,h}(t)
&= h^{q} \frac{{{p-1}\choose{i-1}}}{{{p+q-1}\choose{i+q-1}}}
B_{i+q,p+q,h}(t)
\\
&= h^{q} \frac{{i^{(q)}}}{p^{(q)}}
B_{i+q,p+q,h}(t)
\end{array}
\end{equation}
\end{lemma}
In the notation of the lemma,
\begin{equation} \label{eqdef rising fac}
\begin{array}{r@{}l}
i^{(0)}&=1
\\
i^{(q)}&=i(i+1)(i+2)\cdots(i+q-1)
\end{array}
\end{equation}
is the rising factorial and the following identity was used.
\begin{equation}
\frac{{(i-1)^{(q)}}}{(p-1)^{(q)}}
=
\frac{{{p-1}\choose{i-1}}}{{{p+q-1}\choose{i+q-1}}}
\end{equation}

The following lemma can now be proven.
\begin{lemma} \label{lemBern integral c}
Let $c \geq 0$ and $h >0$. Then, the $L^{2}_{c}$ scalar product of two Bernstein polynomials is given by
\begin{equation}
\begin{array}{r@{}l}
\langle B_{i,p,h}(\. ), B_{j,q,h}(\. ) \rangle_{c}&=
\int _{0}^{h} e^{c\.} B_{i,p,h}(\. ) B_{j,q,h}(\. ) =
\\
&=
{{p-1}\choose{i-1}}{{q-1}\choose{j-1}}
\sum\limits _{n=0}^{\infty}
\frac{c^{n} h^{n+1}}{n!}
\frac{1}{(p+q+n-1){{p+q+n-2}\choose{i+j+n-2}}} 
\\
&=
h \frac{{{p-1}\choose{i-1}}{{q-1}\choose{j-1}}}{(p + q-1) {{p+q-2}\choose{i+j-2}}}
{}_{1}F_{1}(i + j-1, p + q, ch)
\end{array}
\end{equation}
The $L^{2}_{c}$ norm of a Bernstein polynomial is given by
\begin{equation}
\begin{array}{r@{}l}
\| B_{i,p,h}(\. ) \|_{L^{2}_{c}}^{2} &=
\left(
\langle B_{i,p,h}(\. ), B_{i,p,h}(\. ) \rangle_{c}
\right) =
\\
&=
{{p-1}\choose{i-1}}^{2}
\sum\limits _{n=0}^{\infty}
\frac{c^{n} h^{n+1}}{n!}
\frac{1}{(2p+n-1){{2p+n-2}\choose{2i+n-2}}} 
\\
&=
h \frac{{{p-1}\choose{i-1}}^{2}}{(2p-1) {{2p-2}\choose{2i-2}}}
{}_{1}F_{1}(2i-1, 2p, ch)
\end{array}
\end{equation}
\end{lemma}
In statement of the lemma, Kummer's confluent hypergeometric function is used, which
is defined as
\begin{equation}
\begin{array}{r@{}l}
{}_1F_1(a;b;z) &= \sum_{n=0}^\infty \frac {a^{(n)} z^n} {b^{(n)} n!}
\\
&=
\frac{\Gamma(b)}{\Gamma(a)\Gamma(b-a)}\int_0^1 e^{zu}u^{a-1}(1-u)^{b-a-1}\dd u
\end{array}
\end{equation}
where the rising factorial is defined in eq. \eqref{eqdef rising fac}, and the integral
representation holds in the domain $\mathrm{Re}(a)>0$ and $\mathrm{Re}(b)>0$, which is
the case of interest for the lemma.

The proof of lemma uses the power-series expansion of the exponential,
\begin{equation}
e^{t} = \sum _{0}^{\infty}
\frac{1}{n!}t^{n}
\end{equation}
lemma \ref{lemBern prod} and its corollary. Both the lemma and its corollary are, of
course, special cases of this lemma, obtained by posing $c=0$ (where the formal
convention $c^{0} = 1$ when $c=0$ is used).

\begin{proof}
By definition and by lemmas \ref{lemBern integral}, \ref{lemBern prod} and
\ref{lemBer momentum},
\begin{equation}
\begin{array}{r@{}l}
\langle B_{i,p,h}(\. ), B_{j,q,h}(\. ) \rangle_{c}&=
\int _{0}^{h} e^{ct} B_{i,p,h}(t ) B_{j,q,h}(t )\mathrm{d}t
\\
&=
\int _{0}^{h} \sum\limits _{n=0}^{\infty}
\frac{1}{n!}(ct)^{n}
B_{i,p,h}(t ) B_{j,q,h}(t )\mathrm{d}t
\\
&=
\sum\limits _{n=0}^{\infty}
\frac{c^{n}}{n!}
\int _{0}^{h}t^{n}
B_{i,p,h}(t ) B_{j,q,h}(t )\mathrm{d}t
\\
&=
\frac{{{p-1}\choose{i-1}}{{q-1}\choose{j-1}}}{{{p+q-2}\choose{i+j-2}}}
\sum\limits _{n=0}^{\infty}
\frac{c^{n}}{n!}
\int _{0}^{h}
t^{n}
B_{i+j-1,p+q-1,h}(t)
\mathrm{d}t
\\
&=
\frac{{{p-1}\choose{i-1}}
{{q-1}\choose{j-1}}}
{{{p+q-2}\choose{i+j-2}}}
\sum\limits _{n=0}^{\infty}
\frac{c^{n} h^{n}}{n!}
\frac{{{p+q-2}\choose{i+j-2}}}
{{{p+q+n-2}\choose{i+j+n-2}}}
\int _{0}^{h}
B_{i+j+n-1,p+q+n-1,h}(t)
\mathrm{d}t
\\
&=
{{p-1}\choose{i-1}}
{{q-1}\choose{j-1}}
\sum\limits _{n=0}^{\infty}
\frac{c^{n} h^{n+1}}{n!}
\frac{1}{(p+q+n-1){{p+q+n-2}\choose{i+j+n-2}}} 
\end{array}
\end{equation}
The expression using Kummer's confluent hypergeometric function follows from direct
application of the power series definition, or the integral formula and the definition
of the Bernstein polynomials. The series converges as it is bounded above by the
exponential series.
\end{proof}

The proof of the previous lemma also provides a calculation for the product
of an exponential with a Bernstein polynomial and its expression in the
Bernstein polynomial basis. It is stated in the form of a lemma for convenience.

\begin{lemma} \label{lemBern times exp}
    The following formula holds true
    \begin{equation}
    \begin{array}{r@{}l}
        e^{c\.} B_{i,p,h}(\. ) &=
        {{p-1}\choose{i-1}}
        \sum\limits _{n=0}^{\infty}
        \frac{c^{n} h^{n}}{n!}
        \frac{1}
        {{{p+n-1}\choose{i+n-1}}}
        B_{i+n,p+n,h}(t)
        \\
        &=
        \sum\limits _{n=0}^{\infty}
        \frac{c^{n} h^{n}}{n!}
        \frac{i^{(n)}}
        {p^{(n)}}
        B_{i+n,p+n,h}(t)
    \end{array}
    \end{equation}
\end{lemma}
It is reminded that the rising factorial is defined in eq. \eqref{eqdef rising fac}.
The following corollary also follows directly.
\begin{corollary} \label{corBern integral c} 
The integral of a Bernstein polynomial multiplied by an exponential is given by
\begin{equation}
\int _{0} ^{t} e^{c\t} B_{i,p,h}(\t )\dd \t =
\frac{h}{p}
\sum\limits _{n=0}^{\infty}
        \frac{c^{n} h^{n}}{n!}
        \frac{i^{(n)}}
        {(p+1)^{(n)}}
        \sum _{j=i+n+1}^{p+n+1} B_{j,p+n+1,h}(t)
\end{equation}
and its definite integral by
\begin{equation}
\begin{array}{r@{}l}
\int _{0} ^{h}\int _{0} ^{t}  e^{c\t} B_{i,p,h}(\t )\dd \t \dd t &=
\frac{h^{2}(p-i+1)}{p(p+1)}
\sum\limits _{n=0}^{\infty}
        \frac{c^{n} h^{n}}{n!}
        \frac{i^{(n)}}
        {(p+2)^{(n)}}
\\
&=
\frac{h^{2}(p-i+1)}{p(p+1)} {}_{1}F_{1}(i, p+2, ch)
\end{array}
\end{equation}
\end{corollary}

\subsection{Polynomial approximation of functions in Sobolev spaces}

In the following, theorems $2.3$ and $2.4$ 
from \cite{ApproximationLegSob} will be needed. For completeness, the statements and the
context of the needed results are presented here below.

In the proof, the Legendre basis for the space of polynomials of degree $p-1$ is more
useful than the Bernstein polynomials, due to their property of being orthogonal with
respect to the $L^{2}$-norm when restricted in the interval $(-1,1)$. The Legendre
polynomials have a number of possible expressions, both implicit and explicit, one of them
being
\begin{equation}
    L_{n}(x)=
\frac{1}{2^{n} n!}
\frac{\mathrm{d}^{n}}{\mathrm{d}x}
(x^{2}-1)^{n}
\end{equation}
The polynomial $L_{n}(\. )$ is of degree $n$ and satisfies
\begin{equation}
    \begin{array}{r@{}l}
L_{n}(1) &= 1 \\
\int _{-1}^{1} L_{n}(\. ) L_{m}(\. ) &= 0, \forall m \neq n
\\
\| L_{n} \|_{L^{2}}^{2} &= \frac{2}{2n+1}
\end{array}
\end{equation}
so that
\begin{equation} \label{eqBern normalized}
    \phi _{n}(\. ) = \sqrt{\frac{2n+1}{2}}L_{n}(\. ), n \in \N
\end{equation}
form a complete orthonormal basis for $L^{2}(-1,1)$.

For $p\geq 1$, the Legendre basis is given by
\begin{equation}
    \{ \phi _{n} \}_{n = 0, \cdots p-1}
\end{equation}
and is thus equivalent to the Bernstein basis, obtained by
$\{ B_{i,p,1}(\.) \}_{i=1,\cdots,p}$ after a transformation of the domain of the
$B_{i,p,1}$ (which are defined on $[0,1]$), see \cite{FAROUKI2000145} for a study of
the basis transformations. 

A function $u \in L^{2}(-1,1)$ admits the Fourier-Legendre representation
\begin{equation}
    u(\. ) = \sum _{0}^{\infty} u_{n}\phi _{n}(\. )
\end{equation}
where $u_{k} = \int _{-1}^{1}u(\. )\phi _{n}(\. )$ and the orthogonal projection on
the space of polynomials of degree $p-1$ is given by
\begin{equation}
    T_{p}u(\. ) = \sum _{0}^{p-1} u_{n}\phi _{n}(\. )
\end{equation}

With this notation, the following statements hold. They consist in a simple rephrasing of
the corresponding statements in the reference \cite{ApproximationLegSob}. To each statement
a corollary is given, obtained just by transforming the interval $[-1,1]$ to the interval
$[0,h]$, which is of interest in the present article. The proof of the corollaries is
omitted, as it follows directly from the transformation of functions and integrals under
composition with affine transformations.
\begin{theorem}[Theorem $2.3$ of \cite{ApproximationLegSob}]
For any real $\s \geq 0$, there exists a constant $C$ such that for any $p \in \N $,
\begin{equation}
    \| u - T_{p}u\|_{L^{2}(-1,1)} \leq C (p+1)^{-\s}  \| u\|_{H^{\s}(-1,1)},
\forall u \in H^{\s}(-1,1)
\end{equation}
\end{theorem}

\begin{corollary} \label{corproj error}
For any real $\s \geq 0$, there exists a constant $C$ such that for any $p \in \N $,
\begin{equation}
    \| u - T_{p}u\|_{L^{2}(0,h)} \leq C (p+1)^{-\s} h^{\s} \| u\|_{H^{\s}(0,h)}
\end{equation}
for all $u \in H^{\s}(0,h)$ such that $\int _{0}^{h}u(\.) = 0$.
\end{corollary}


\begin{theorem}[Theorem $2.4$ of \cite{ApproximationLegSob}] \label{propApproxSob}
For any real $0 \leq \m \leq \s$, there exists a constant $C$ such that for any $p \in \N $,
\begin{equation}
    \| u - T_{p}u\|_{H^{\m}(-1,1)} \leq C (p+1)^{e(\m, \s)}  \| u\|_{H^{\s}(-1,1)},
\forall u \in H^{\s}(-1,1)
\end{equation}
where
\begin{equation}
    e(\m, \s) = \begin{cases}
2\m - \s - 1/2, \m \geq 1
\\
3\m /2 - \s, 0 \leq \m \leq 1
\end{cases}
\end{equation}
\end{theorem}

\begin{corollary} \label{corproj error interpol}
For any real $\s \geq 0$, there exists a constant $C$ such that for any $p \in \N $,
\begin{equation}
    \| u - T_{p}u\|_{H^{\m}(0,h)} \leq C (p+1)^{e(\m, \s)}h^{\s-\m}  \| u\|_{H^{\s}(0,h)},
\forall u \in H^{\s}(0,h)
\end{equation}
for all $u \in H^{\s}(0,h)$ such that $\int _{0}^{h}u(\.) = 0$.
\end{corollary}

Finally, the following lemma (lem. $2.4$ of
\cite{ApproximationLegSob}) on the growth of Sobolev
norms of polynomials and its corollary will be needed.

\begin{lemma}
    For $p \in \N^{*}$ and $u$ a polynomial of degree
    $p$, i.e. such that $T_{p}u = u$, the following
    holds
    \begin{equation}
        \|u\|_{H^{\m}(-1,-1)} \leq C p ^{2(\m - \nu)}
        \|u\|_{H^{\n}(-1,-1)}
    \end{equation}
\end{lemma}

\begin{corollary}
    For $p \in \N^{*}$, $h>0$ and $u$ a polynomial of degree
    $p$, i.e. such that $T_{p}u = u$, the following
    holds
    \begin{equation}
        \|u\|_{H^{\m}(0,h)} \leq C p ^{2(\m - \nu)}
        h^{\m - \n}
        \|u\|_{H^{\n}(0,h)}
    \end{equation}
\end{corollary}


\section{From Legendre to Bernstein basis}

In this section, the Legendre polynomials defined on the interval $[0,1]$ will be
used, given by
\begin{equation}
\tilde{L}_{n}(x) = L_{n}(2x-1)
\end{equation}
The following formula from \cite{FAROUKI2000145} is a restatement of the
proposition of the paper using the notation established herein.

\begin{prop}[\cite{FAROUKI2000145}, Prop. $2$] \label{propBernLeg}
Let
\begin{equation}
B_{i,p,1}(\. ) = \sum _{m=0}^{p-1} \sqrt{2m+1} \Lambda ^{-1} _{im}(p)\tilde{L}_{m}(\.)
\end{equation}
be the expression of a Bernstein polynomial as a linear combination of Lagrange
polynomials, where the $\sqrt{2m+1}$ factor imposes
$L^{2}$ normalization. Then
\begin{equation}
\Lambda ^{-1} _{im} = \Lambda ^{-1} _{im} (p) = \frac{\sqrt{2m+1}}{p+m}{{p-1}\choose{i-1}}
\sum _{q=0}^{m}(-1)^{m+q}\frac{{{m}\choose{q}}^{2}}{{{p+m-1}\choose{i+q-1}}}
\end{equation}
where $1 \leq i \leq p$ and $0 \leq m \leq p-1$.
\end{prop}

Using the expression obtained in lem. \ref{lemBern times exp}, the following corollary is
immediate.
\begin{corollary}
For a Bernstein polynomial multiplied by an exponential, it holds that
\begin{equation}
e^{c\.}B_{i,p,1}(\. )
=
\sum _{m=0}^{\infty} \sqrt{2m+1} \mathcal{M} ^{-1} _{im}(p)\tilde{L}_{m}(\.)
\end{equation}
where
\begin{equation}
\begin{array}{r@{}l}
\mathcal{M} ^{-1} _{im}(p) &=
\sum _{n\geq \max \{m-p-1,0 \}}
\frac{c^{n}h^{n}}{n!} 
\frac{{{p-1}\choose{i-1}}}{{{p+n-1}\choose{i+n-1}}}
\Lambda ^{-1} _{(i+n),m} (p+n)
\\
&=
\sum _{n\geq \max \{m-p-1,0 \}}
\frac{c^{n}h^{n}}{n!}  \frac{{{p-1}\choose{i-1}}}{p+m+n}
\sum _{q=0}^{m} (-1)^{m+q}\frac{{{m}\choose{q}}^{2}}{{{p+n+m-1}\choose{i+n+q-1}}}
\end{array}
\end{equation}
\end{corollary}

Anticipating \S \ref{secerror estimates}, a brief study on the higher-order
truncation of Bernstein polynomials
in the Legendre basis will be will be presented below. In \cite{FAROUKI2000145},
the author discusses the difficulty of simplifying formulas like the
one in prop. \ref{propBernLeg} here above. The authors of the present
work tried to obtain at least estimates for the decay of
$\max _{i} |\Lambda ^{-1} _{im}|$ for $m \in \{ p-2, p-1 \}$, which
would be useful in the error estimates, but no exploitable result was obtained.
What was opted for instead was a numerical study covering the range of
useful degrees of approximation $p$, which does not exceed $25$ due to
numerical stability issues.

\begin{figure}[h!]
  \centering
	  \includegraphics[width=0.7\linewidth]{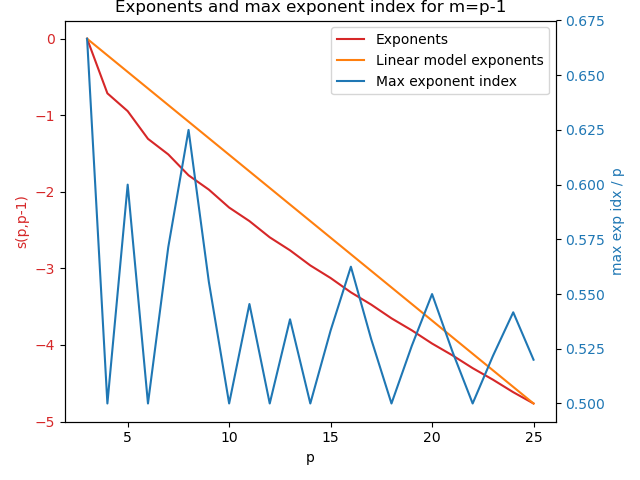}
  \caption{The exponent expressing the size of the highest order coefficient of a
  Bernstein polynomial written in the Legendre basis as a power of $p$.}
  \label{fig:exp_last}
\end{figure}

\begin{figure}[h!]
  \centering
  \includegraphics[width=0.7\linewidth]{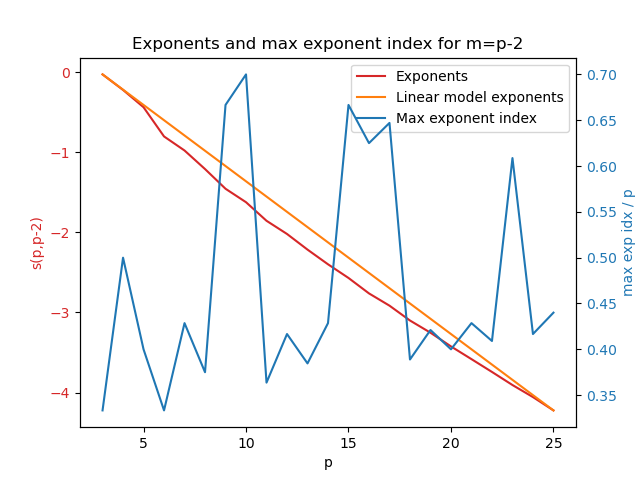}
  \caption{The exponent expressing the size of the second highest order coefficient of
  a Bernstein polynomial written in the Legendre basis as a power of $m=p$.}
  \label{fig:exp_next_to_last}
\end{figure}

\begin{figure}[h!]
  \centering
  \includegraphics[width=0.7\linewidth]{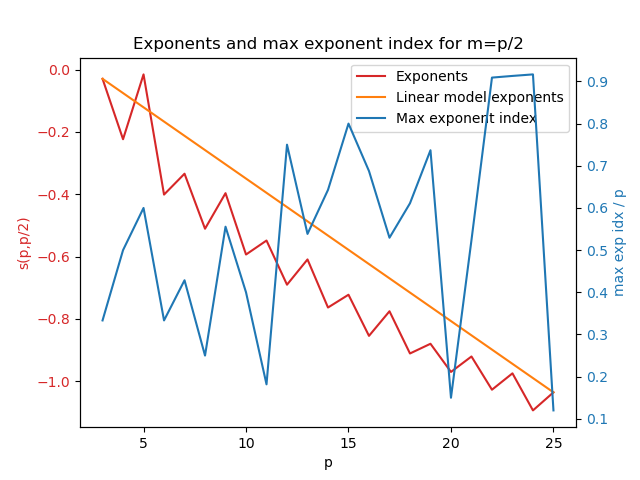}
  \caption{The exponent expressing the size of a middle order coefficient of a
  Bernstein polynomial written in the Legendre basis as a power of
  $m=p$.}
  \label{fig:exp_low}
\end{figure}

The following graphs in figures \ref{fig:exp_last}, \ref{fig:exp_next_to_last} and
\ref{fig:exp_low} feature the exponent $s$ satisfying
\begin{equation}
    s = s(p,m) = \max _{i} \frac{\log  |\Lambda ^{-1} _{im}| - \log \| B_{i,p,1} - \frac{1}{p}\|_{0}
    }{\log p}
\end{equation}
so that
\begin{equation}
    |\Lambda ^{-1} _{im}| \leq
    \| B_{i,p,1} - \frac{1}{p}\|_{0} p^{s},
    \forall 1 \leq i \leq p
\end{equation}
for $m = p-2,p-1$, where the homogeneous $L^{2}$ norm
$\| B_{i,p,1} - \frac{1}{p}\|_{0}$ is obtained by lem.
\ref{lemBern prod} using
\begin{equation}
    \| B_{i,p,1} - \frac{1}{p}\|_{0}^{2} = \| B_{i,p,1}\|_{0}^{2} - \frac{1}{p^{2}}
\end{equation}
The growth is compared with $p$ even though the degree of
$B_{i,p,1}$ is $p-1$ and therefore its norm growth is
governed by powers of $p-1$, purely for reasons of
convenience.
This quantity will be relevant in the error rate estimates. The
case of a lower degree, $m = \lfloor p/2 \rfloor $ is also presented for
reasons of comparison, in order to show that lower degree Legendre polynomials do
bear a larger part of the $L^{2}$ norm of the Bernstein polynomials, as is
expected from general theory.

A linear model giving an upper bound for $s$ as a function of $m$ is also
given, for illustration purposes. The linear models, fitting the actual data
rather tightly, read
\begin{equation}
    \begin{array}{r@{}l}
         s(p,p-1) &= -0.216573775474902 p + 0.649721326424705,
         \\
         s(p,p-2) &= -0.190609494287812 p + 0.542455736079486,
         \\
         s(p,\lfloor p/2 \rfloor) &= -0.0457451923495941 p + 0.107862830264831  \\
    \end{array}
\end{equation}
It is clear from the graphs, as well as from the linear upper bounds, that the
lower-order coefficient bears a more significant part of the weight than the higher order
coefficients, and that its dependence on $p$ is milder,
as should be expected from theory.

\section{From $[0,h]$ to $[0,\T]$, the space $\PP_{h,p}$}

All the above properties are proved for the basis defined on the interval $[0,h]$. Since
the goal is to obtain an approximate solution defined on $[0,\T]$, one can translate the
functions $B_{i}(\. )$ by $j\. h$ for $0 \leq j \leq l-1$, in order to obtain the functions
\begin{equation}
B_{i}^{j}(t) = B_{i}(t-jh) = \begin{cases}
\frac{\Gamma (p)}{\Gamma (i) \Gamma (p-i +1)}
\left( \frac{t - jh}{h}
\right)^{i-1}
\left( \frac{h + jh-t}{h}
\right)^{p-i} \text{ if } t \in I_{j}\\
0 \hphantom{
\frac{\Gamma (p)}{\Gamma (i) -}
\left( \frac{t}{h}
\right)^{i-1}
} \text{   otherwise}
\end{cases}
\end{equation}
where $B_{i}^{0}(t) \equiv B_{i}(t)$. The functions
\begin{equation}
\{ B_{i}^{j}(\. ) , 1 \leq i \leq p , 0 \leq j \leq l-1 \}
\end{equation}
form a basis of the piecewise polynomial functions of degree $p$ on $\bar{I}$, with
break-points $\{ j h \}_{j=1}^{l-1}$. In other words, if $\PP _{h,p}$ is the space of
such functions, $f \in \PP _{h,p}$ if, and only if, there exist $f_{i}^{j} \in \R$ such
that
\begin{equation} \label{eq expr basis non hom}
f(\. ) = \sum _{i,j} f_{i}^{j} B_{i}^{j}(\. )
\end{equation}
It should be noted that the functions in $\PP _{h,p}$ are not well defined at break-points, but
since the break-points form a set of measure $0$ this is not of concern. The ambiguity
in the definition will be lifted rightaway when continuity conditions at break-points
are imposed. The dimension of the space $\PP _{h,p}$ is equal to $l\.p$.

Piecewise polynomial functions are in $H^{1}$ iff they are in $C^{1}$. The SDOF equation
being a second order equation, the value of the first derivative is an independent
variable. As a consequence, it has to be determined by the choice of functions
at break-points and not by the algorithm itself. The space in which the solution should be
placed is, therefore,
\begin{equation}
\PP _{h,p}^{1} = \PP _{h,p} \bigcap C^{1}(I)
\end{equation}
A function $f \in \PP _{h,p}$ is $C^{1}$ iff the function itself and its
first derivative are continuous at break-points. In terms of the coefficients $f_{i}^{j}$,
this is equivalent to
\begin{equation} \label{eq cond cont}
\begin{array}{c@{}l}
f_{1}^{j+1}&= f_{p}^{j} \\
f_{2}^{j+1} &= - f_{p-1}^{j} +2f_{p}^{j}
\end{array}
\end{equation}
The space $\PP _{h,p}^{1} $ is consequently of dimension $l\. p - 2(l-1) $, since it is obtained
by imposing $2$ linearly independent conditions on each break-point other than $0$ and $\T$
on functions in $\PP _{h,p} $.

Since the acceleration is a dependent variable, no $C^{2}$ continuity conditions are
imposed at break-points. The authors actually tried the version of the algorithm where such
conditions were imposed, and the method failed to converge.

%

\section{The first step of the algorithm, $j=0$} \label{secfirst step}

In what follows, $p\geq 3$ is an arbitrary natural number, and $l=1$ so that $h=\T$.

The approximate solution $x_{ap}$ restricted to the interval $[0,h]$ is thus assumed to have the form
\begin{equation}
x_{ap} (\. )= \sum _{1\leq i \leq p} u^{0}_{i} B^{0}_{i}(\. )
= \sum _{1\leq i \leq p} u_{i} B_{i}(\. )
\end{equation}
where $u_{i} \in \R$ are the $p$ unknown coefficients. The function $v_{0,h}$ bearing the boundary conditions
is defined as
\begin{equation}
v_{0,h} = x_{0}B_{1}(\. ) + x_{h }B_{p}(\. )
\end{equation}

Since only the first and last functions bear the boundary conditions, we immediately obtain
\begin{equation}
\begin{cases}
u_{1} = x_{0} \\
u_{p}  = x_{h }
\end{cases}
\end{equation}
where $x_{h }$ is unknown and treated as a parameter.

Consequently, the approximate solution to the weak problem, $u_{ap}(\. )$ assumes the form
\begin{equation}
u_{ap} (\. )= \sum _{1< i < p} u_{i} B_{i}(\. )
\end{equation}
so that
\begin{equation}
x_{ap} (\. )= u _{ap} (\. ) + v_{0,h} (\. )
\end{equation}

Direct application of the method outlined in \S 9 of \cite{NKSDOFI} gives the
$(p-2)\times p$ matrix $[\B ]$ with elements
\begin{equation} \label{eqmatrix B}
\B _{ij } = - \langle
\dot{B}_{i}(\. ) , \dot{B}_{j }(\. )
\rangle _{c}
+ k \langle
B_{i}(\. ) , B_{j }(\. )
\rangle _{c}
, 2 \leq i  \leq p-1
, 1 \leq j  \leq p
\end{equation}
and the $(p-2)$ column vector $\F $ with elements
\begin{equation} \label{eq def vec F}
\F _{i} = \langle
f(\. ) , B_{i}(\. )
\rangle _{c} -
 \langle
u_{0} \dot{B}_{1}(\. ) + u_{p } \dot{B}_{p}(\. ) , \dot{B}_{i}(\. )
\rangle  _{c} +
k \langle
u_{0}B_{1}(\. ) + u_{p } B_{p}(\. ) , B_{i}(\. )
\rangle _{c}
\end{equation}
for $2 \leq i  \leq p-1$, depending on the parameter $u_{p}$.

The linear system
\begin{equation} \label{eq lin sys simple}
[\tilde{\B} ].\{ u_{i} \}_{i=2}^{p-1}=\{ \F \}
\end{equation}
where $\tilde{\B}$ is obtained by dropping the first and last columns of matrix
$\B$, solves the problem assuming known boundary conditions and thus treating
$u_{p}$ as a parameter. Then, the equation
\begin{equation} \label{eqinitial velocity Bernstein}
\dot{u}(0)=\dot{x}_{0} \iff  u_{1} = \frac{h}{p-1}\dot{x}_{0} + x_{0}
\end{equation}
allows to solve for $u_{p}$ and determine the rest of the coefficients $u_{i}$.

In total, if one now treats $u_{1}$ and $u_{p}$ as unknown quantities and moves the terms in $\{ \F \}$
containing them to the lhs of eq. \eqref{eq lin sys simple} and considers the full vector of coefficients
$\{ u_{i} \}_{i=1}^{p}$, they obtain the following linear system
\begin{equation} \label{eq lin sys simple aug}
[\B \B ].\{ u_{i} \}_{i=1}^{p}=\{ \F \F \}
\end{equation}
where the matrix of known coefficients
\begin{equation}
[\B \B ]_{p\times p}
=
\left[
\begin{array}{c|c}
\begin{bmatrix}
1 & 0 \\
-\frac{p-1}{h} & \frac{p-1}{h}
\end{bmatrix}
& [ 0 ]_{2\times(p-2)} \\
\midrule
\{ \B ^{(1)} \}_{(p-2)\times 2} &
[\B ^{(2)} ]_{(p-2)\times (p-2)} 
\end{array}
\right]
\end{equation}
with $ \B ^{(1)}$ formed by the first two columns of matrix $\B$, and $ \B ^{(2)}$
by the remaining $p-2$ columns of the same matrix. The rhs vector reads
\begin{equation} \label{eq vec FF}
\{ \F \F \}_{(p)\times 1}
=
\begin{Bmatrix}
x_{0} \\
\dot{x}_{0} \\
\{ \langle
f(\. ) , B_{i}(\. )
\rangle _{c} \}_{i=2}^{p-1} \}
\end{Bmatrix}
\end{equation}
or, in a more compact form,
\begin{equation} \label{eq BB}
[\B \B ]_{m\times p}
=
\left[
\begin{array}{c|c}
[C]_{2\times 2}
& [ 0 ]_{2\times(p-2)} \\
\midrule
\{ \B ^{(1)} \}_{(p-2)\times 2} &
[\B^{(2)} ]_{(p-2)\times (p-2)} 
\end{array}
\right]
\end{equation}
where
\begin{equation}
[C]_{2\times 2} =
\begin{bmatrix}
1 & 0 \\
-\frac{p-1}{h} & \frac{p-1}{h}
\end{bmatrix}
\end{equation}
%
%
%

The linear system can be efficiently solved in the following way:
\begin{enumerate}
\item the first two equations are decoupled from the rest of the system and can
be solved separately. They are already in (lower) triangular form.
\item the first column of the vector $\{ \B ^{(1)} \}_{(p-2)\times 1}$ multiplied by $u_{1}= x_{0}$ and
its second column multiplied by $u_{2}$ are subtracted from the vector
$\{ \langle f(\. ) , B_{i}(\. ) \rangle _{c} \}_{i=2}^{p-1} \}$ giving the vector
$\{ \tilde{\F}  \} $
\item the rest of the coefficients, $\{ u_{i} \}_{i=3}^{p}$, are obtained by
\begin{equation}
\{ u_{i} \}_{i=3}^{p} = [\B^{(2)} ]^{-1}. \{  \tilde{\F}\}
\end{equation}
\end{enumerate}

The first step of the algorithm is completed.

\section{The remaining steps}

One then can calculate the final displacement
\begin{equation} \label{eqfinal disp}
x_{ap}(h ) = u_{p}
\end{equation}
and final velocity
\begin{equation} \label{eqfinal speed}
\dot{x}_{ap}(h) = \frac{p-1}{h}(u_{p}-u_{p-1})
\end{equation}
and iterate the algorithm for the desired number of timesteps in the following way.

For the second step, $j = 1$, the approximate solution $x_{ap}$ restricted to the
interval $[h,2h]$ and is thus assumed to have the form
\begin{equation}
x_{ap} (\. )= \sum _{1\leq i \leq p} u^{1}_{i} B^{1}_{i}(\. )
= \sum _{1\leq i \leq p} u_{i}^{1} B_{i}(\. -h )
\end{equation}
The initial displacement and velocity are $x_{ap}(h ) $ and $ \dot{x}_{ap}(h) $ as
calculated here above and the algorithm of the first step applies verbatim. The next
steps are carried out in the same fashion.

\begin{remark}
A direct way to parallelize the algorithm using two processors would be the following. The
first processor calculates and inverts the matrix $[\B ]$, obtaining the matrix
$[\B ]^{-1}$ that is used throughout the algorithm. The second processor calculates the
rhs vector $\{ \F \}$, appending the vector corresponding to the timestep $j+1$ at the end
of the vector constructed at the end of the $j$-th timestep, and feeds the result to the
first processor, which applies the algorithm for each timestep.

There does not seem to be a way to use more processors for this problem. However, piecemeal
computation by the second processor of the vector $\{ \F \}$ for each timestep $j$, instead
of construction of the entire vector of dimension $j*p$ is a direct improvement in terms of
memory usage, and should become necessary when treating large-scale MDOF problems.
\end{remark}

\section{The case $p=3$}

The case where $p=3$ is the only case where a direct comparison can be made between
the method proposed herein with the traditional step-wise methods, as it is the case
of an approximation by a quadradic polynomial. In this case, the matrix
$[\B ]_{(p-2)\times 3}$ of eq. \eqref{eqmatrix B} is a $1\times 3$ vector, and all
calculations can be made explicit.

\subsection{The undamped case}

The expressions when $c>0$, even though still explicit, become cumbersome. The case
$c=0$ merits, thus, a separate study for clarity of exposition.


The matrix of eq. \eqref{eq lin sys simple aug} reads
\begin{equation}
[\B \B ]_{3\times 3}
=
\left[
\begin{array}{c|c}
\begin{bmatrix}
1 & 0 \\
-\frac{2}{h} & \frac{2}{h}
\end{bmatrix}
& \begin{bmatrix}
0 \\
0
\end{bmatrix} \\
\midrule
\{ \frac{2}{3h} + \frac{kh}{10},
-\frac{4}{3h} + \frac{2kh}{15} \} &
[\frac{2}{3h}+\frac{kh}{10}]
\end{array}
\right]
\end{equation}
or, equivalently, already in lower triangular form
\begin{equation}
[\B \B ]_{3\times 3}
=
\begin{bmatrix}
1 & 0 & 0 \\
-\frac{2}{h} & \frac{2}{h} & 0 
\\
\frac{2}{3h} + \frac{kh}{10}
&
-\frac{4}{3h} + \frac{2kh}{15}
&
\frac{2}{3h}+\frac{kh}{10}
\end{bmatrix}
\end{equation}
The rhs vector of the linear system as defined in eq. \eqref{eq vec FF} reads
\begin{equation}
\{ \F \F \}_{3\times 1}
=
\begin{Bmatrix}
x_{0} \\
\dot{x}_{0} \\
\langle
f(\. ) , B_{2.3,h}(\. )
\rangle _{0} 
\end{Bmatrix}
\end{equation}

Consequently, the linear system $[\B \B ] \{ u_{i}\}_{1}^{3} = \{ \F\F\}$ admits the
unique solution
\begin{equation}
\begin{array}{r@{}l}
u_{1} &= x_{0}
\\
u_{2} &= \frac{h}{2}\dot{x}_{0} + x_{0}
\end{array}
\end{equation}
which yields
\begin{equation}
u_{3}
= \frac{30h\langle
f(\. ) , B_{2,3,h}(\. )
\rangle _{0}
+
(20 - 7kh^{2})x_{0}
+
2(10 - kh^{2})h\dx_{0}
}
{20+3kh^{2}}
\end{equation}

This results in the final displacement at time $t=h$, $x_{h}$ being equal to $u_{3}$
as implied by eq. \eqref{eqfinal disp},
i.e. by
\begin{equation}
x_{h} = 
\frac{30h\langle
f(\. ) , B_{2,3,h}(\. )
\rangle _{0}
+
(20 - 7kh^{2})x_{0}
+
2(10 - kh^{2})h\dx_{0}
}
{20+3kh^{2}}
\end{equation}
The final velocity, $\dx_{h}$, as implied by eq. \eqref{eqfinal speed}, is given by
\begin{equation} \label{eqvelocity quadratic}
\dx_{h} = \frac{60\langle
f(\. ) , B_{2,3,h}(\. )
\rangle _{0}
-
20khx_{0}
+
(20 - 7kh^{2})\dx_{0}}{20+3kh^{2}}
\end{equation}

Using these formulas, the following theorem can be proved.
\begin{theorem}
For $c=0$, $k>0$, and $p=3$, the algorithm proposed in \S \ref{secfirst step} boils down
to the step-wise iteration of
\begin{equation}
\begin{cases}
x_{j+1} = 
\frac{30h\int _{jh}^{(j+1)h}
f(\. ) B_{2,3,h}(\. )
+
(20 - 7kh^{2})x_{j}
+
2(10 - kh^{2})h\dx_{j}
}
{20+3kh^{2}}
\\
\dx_{j+1} = \frac{60 \int _{jh}^{(j+1)h}
f(\. ) B_{2,3,h}(\. )
-
20khx_{j}
+
(20 - 7kh^{2})\dx_{j}}{20+3kh^{2}}
\end{cases}
\end{equation}
where $x_{0}$ and $\dx_{0}$ are given as initial conditions, and $f$ is the external
excitation force.

The approximation rates are as follows:
\begin{enumerate}
\item For $f \equiv 0$,
\begin{enumerate}
\item for $x_{j}=1$ and $\dx _{j}=0$, 
\begin{equation}
\begin{array}{r@{}l}
|x_{j+1}-x_{ex}((j+1)h)| &= O(h^{4})
\\
|\dx_{j+1}-\dx_{ex}((j+1)h)| &= O(h^{3})
\\
|\mathrm{ME}_{j+1}-\mathrm{ME}_{ex}((j+1)h)| &= O(h^{4})
\end{array}
\end{equation}
\item for $x_{j}=0$ and $\dx _{j}=1$, 
\begin{equation}
\begin{array}{r@{}l}
|x_{j+1}-x_{ex}((j+1)h)| &= O(h^{3})
\\
|\dx_{j+1}-\dx_{ex}((j+1)h)| &= O(h^{4})
\\
|\mathrm{ME}_{j+1}-\mathrm{ME}_{ex}((j+1)h)| &= O(h^{4})
\end{array}
\end{equation}
\end{enumerate}
\item For $x_{0}=\dx_{0}=0$ and $f$ piece-wise constant in each interval $[jh,(j+1)h]$:
\begin{equation}
\begin{array}{r@{}l}
|x_{j+1}-x_{ex}((j+1)h)| &= O(h^{4})
\\
|\dx_{j+1}-\dx_{ex}((j+1)h)| &= O(h^{3})
\\
|\mathrm{ME}_{j+1}-\mathrm{ME}_{ex}((j+1)h)| &= O(h^{4})
\end{array}
\end{equation}
The constant in the big-O notation is proportional to $f_{j} = f \big|_{[jh,(j+1)h]}$.
\end{enumerate}
\end{theorem}
In the notation of the theorem, $x_{ex}$ represents the exact solution of the
corresponding SDOF problem,
\begin{equation}
\mathrm{ME}_{j} = \frac{1}{2}kx_{j}^{2} + \frac{1}{2}\dx_{j}^{2}
\end{equation}
is the mechanical energy of the approximate solution, and $\mathrm{ME}_{ex}(jh)$
is the mechanical energy of the exact solution at time $t=jh$.

\begin{proof}
The iteration is just a restatement of the calculations preceding the statement of the
theorem.

Regarding the error rates, they follow from the calculation of the Taylor expansions 
of the iteration formulas for the corresponding cases.

More precisely, the Taylor expansion of the function
\begin{equation}
\frac{20 - 7t^{2}}{20+3t^{2}} - \cos t
\end{equation}
around $t=0$ is
\begin{equation} \label{eqexpansion cos}
\frac{t^{4}}{30} - \frac{71t^{6}}{7200} + O(t^{8})
\end{equation}
Similarly, the expansion of
\begin{equation}
\frac{2(10 - t^{2})t}{20+3t^{2}} - \sin t
\end{equation}
is
\begin{equation}
-\frac{t^{3}}{12} - \frac{7t^{5}}{240} + O(t^{7})
\end{equation}
These Taylor expansions account for the good approximation of the displacement in
a homogeneous problem, with $0$ excitation force and non-trivial initial conditions
(displacement and velocity, respectively).

In the same fashion, the expansion of
\begin{equation} \label{eqtaylor sine}
\frac{-20t}{20+3t^{2}} + \sin t
\end{equation}
reads
\begin{equation}
-\frac{t^{3}}{60} - \frac{17t^{5}}{200} + O(t^{7})
\end{equation}
which, together with the expansion of eq. \eqref{eqexpansion cos}, accounts for
the good approximation of the velocity in the same problem.

Regarding the conservation of Mechanical Energy, the Taylor expansion of
\begin{equation}
\left(
\frac{20 - 7t^{2}}{20+3t^{2}}
\right)^{2}
+
\left(
\frac{20t}{20+3t^{2}}
\right)^{2} - 1
\end{equation}
is
\begin{equation}
\frac{t^{4}}{40} - \frac{3t^{6}}{100} + O(t^{8})
\end{equation}
while that of
\begin{equation}
\left(
\frac{2(10 - t^{2})t}{20+3t^{2}}
\right)^{2}
+
\left(
\frac{20 - 7t^{2}}{20+3t^{2}}
\right)^{2}
- 1
\end{equation}
is
\begin{equation}
-\frac{t^{4}}{10} + \frac{t^{6}}{25} + O(t^{8})
\end{equation}
These approximations establish the good energy-conservation properties of the
simplest case of the algorithm proposed herein.

Concerning the terms containing the external force, the following can be obtained
by direct calculation. One gets, for $F(t) = \int _{0}^{t} f(t)  \dd \t$,
\begin{equation}
\begin{array}{r@{}l}
\langle
f(\. ) , B_{2,3,h}(\. )
\rangle _{0}
&=
\int _{0}^{h}
f(t)B_{2,3,h}(t) \dd t
\\
&=
-\frac{2}{h}
\int _{0}^{h}
F(t)\dot{B}_{2,3,h}(t) \dd t
\\
&=
-\frac{2}{h}
\int _{0}^{h}
F(t)(B_{1,2,h}(t) - B_{2,2,h}(t)) \dd t
\\
&=
\frac{2}{h^{2}}
\int _{0}^{h}
F(t)(2t-1) \dd t
\end{array}
\end{equation}
Since $2t-1$ is the Legendre polynomial of order $1$, this last expression is  the 
projection of $F$ onto linear polynomials, where it is reminded that $F$ is defined
modulo a constant. Equivalently, $\langle f(\. ) , B_{2,3,h}(\. ) \rangle _{0}$ is the
projection of $f \in H^{-1}$ onto constants. Replacing $f(\.) \equiv 1$ in
$\langle f(\. ) , B_{2,3,h}(\. ) \rangle _{0}$ gives
\begin{equation}
\langle
1 , B_{2,3,h}(\. )
\rangle _{0}
=
\int _{0}^{h}
B_{2,3,h}(t) \dd t
= \frac{h}{3}
\end{equation}

The solution of the undamped SDOF system for $k=1$, homogeneous initial conditions and
$f(\.) \equiv 1$ reads
\begin{equation}
x(t) = 1-\cos t
\end{equation}
The Taylor development around $0$ of
\begin{equation}
\frac{10t^{2}}{20+3t^{2}} - (1-\cos t)
\end{equation}
reads
is
\begin{equation}
-\frac{t^{4}}{30} + \frac{71t^{6}}{7200} + O(t^{8})
\end{equation}
which establishes that the method produces a $4$th order approximation under the
assumption that the external excitation is constant in $[0,h]$. Concerning the
velocity, the force term in eq. \eqref{eqvelocity quadratic} for $f(\.) \equiv 1$ reads
$20h$. This term gives rise to the same Taylor development as in eq.
\eqref{eqtaylor sine}, resulting in a $3$rd order approximation for the velocity.

The Taylor development for the Mechanical Energy of the system reads
\begin{equation}
\frac{1}{60}h^{4} - \frac{37}{7200}h^{6} +O(h^{8})
\end{equation}
\end{proof}

\subsection{The damped case}


In the case where $c\geq 0$, the matrix of eq. \eqref{eq lin sys simple aug} reads
\begin{equation}
[\B \B ]_{3\times 3}
=
\begin{bmatrix}
1 & 0 & 0 \\
-\frac{2}{h} & \frac{2}{h} & 0 
\\
\Xi _{1}
&
\Xi _{2}
&
\Xi _{3}
\end{bmatrix}
\end{equation}
where
\begin{equation} \label{eqdef Xi}
\begin{array}{r@{}l}
\Xi _{1} =\Xi _{1}(h,c,k) &= \frac{2}{3h}(2{}_{1}F_{1}(1,4,ch)-{}_{1}F_{1}(2,4,ch))
 + \frac{kh}{10}{}_{1}F_{1}(2,6,ch)
\\
\Xi_{2} = \Xi_{2}(h,c,k) &= -\frac{4}{3h}(2{}_{1}F_{1}(1,4,ch)-{}_{1}F_{1}(2,4,ch) + {}_{1}F_{1}(3,4,ch))
\\
&\phantom{\frac{4}{3h}(2{}_{1}F_{1}(1,4,ch)-{}_{1}F_{1}(2,4,ch)}
 + \frac{2kh}{15}{}_{1}F_{1}(3,6,ch)
\\
\Xi_{3} = \Xi_{3}(h,c,k) &= -\frac{2}{3h}({}_{1}F_{1}(2,4,ch) -2 {}_{1}F_{1}(3,4,ch))
 + \frac{kh}{10}{}_{1}F_{1}(4,6,ch) 
\end{array}
\end{equation}
Since ${}_{1}F_{1}(i,p,0) = 1$ for all relevant values of $i,p$, one obtains directly
the formulas for the undamped case when $c$ is set to $0$.

The rhs vector of the linear system as defined in eq. \eqref{eq vec FF} reads
\begin{equation}
\{ \F \F \}_{3\times 1}
=
\begin{Bmatrix}
x_{0} \\
\dot{x}_{0} \\
\int _{0}^{h}
e^{ct} f(t)B_{2.3,h}(t)
\end{Bmatrix}
\end{equation}

This yields the solution
\begin{equation}
\begin{array}{r@{}l}
u_{1} &= x_{0}
\\
u_{2} &= \frac{h}{2}\dot{x}_{0} + x_{0}
\\
u_{3}
&= \frac{
\int _{0}^{h}
e^{ct} f(t)B_{2.3,h}(t)
-
(\Xi_{1}+\Xi_{2})x_{0}
-
\frac{h}{2}\Xi_{2}\dx_{0}
}
{\Xi_{3}}
\end{array}
\end{equation}
Explicit expressions for the quantities $\Xi_{i}$, $i=1,2,3$ can be obtained, but
they are cumbersome and not actually useful.

As in the previous paragraph, the displacement and velocity at time $t=h$ can be
obtained by the following expressions
\begin{equation}
\begin{array}{r@{}l}
x_{h} 
&= \frac{
\int _{0}^{h}
e^{ct} f(t)B_{2.3,h}(t)
-
(\Xi_{1}+\Xi_{2})x_{0}
-
\frac{h}{2}\Xi_{2}\dx_{0}
}
{\Xi_{3}}
\\
\dx_{h}
&=
\frac{2}{h}
\frac{
\int _{0}^{h}
e^{ct} f(t)B_{2.3,h}(t)
-
(\Xi_{1}+\Xi_{2}+\Xi_{3})x_{0}
-
\frac{h}{2}(\Xi_{2}+\Xi_{3})\dx_{0}
}
{\Xi_{3}}
\end{array}
\end{equation}
The following theorem can be now be proved.
\begin{theorem}
For $c\geq 0$, $k>0$, and $p=3$, the algorithm proposed in \S \ref{secfirst step} boils
down to the step-wise iteration of
\begin{equation}
\begin{cases}
x_{j+1} = 
\frac{\int _{jh}^{(j+1)h}
e^{c\.}
f(\. ) , B_{2,3}(\. )
-
(\Xi_{1}+\Xi_{2})x_{j}
-
\frac{h}{2}\Xi_{2}\dx_{j}
}
{\Xi_{3}}
\\
\dx_{j+1} = \frac{2}{h}
\frac{\int _{jh}^{(j+1)h}
f(\. ) , B_{2,3}(\. )
-
(\Xi_{1}+\Xi_{2}+\Xi_{3})x_{0}
-
\frac{h}{2}(\Xi_{2}+\Xi_{3})\dx_{0}
}
{\Xi_{3}}
\end{cases}
\end{equation}
where $x_{0}$ and $\dx_{0}$ are given as initial conditions, and $f$ is the external
excitation force and the $\Xi_{i}$ are defined in eq. \eqref{eqdef Xi}.

The approximation rates are as follows:
\begin{enumerate}
\item For $f \equiv 0$,
\begin{enumerate}
\item for $x_{j}=1$ and $\dx _{j}=0$, 
\begin{equation}
\begin{array}{r@{}l}
|x_{j+1}-x_{ex}((j+1)h)| &= O(ch^{3})
\\
|\dx_{j+1}-\dx_{ex}((j+1)h)| &= O(h^{3})
\\
|\mathrm{ME}_{j+1}-\mathrm{ME}_{ex}((j+1)h)| &= O(ch^{3})
\\
|\mathrm{ME}^{c}_{j+1}-\mathrm{ME}^{c}_{ex}((j+1)h)| &= O(ch^{3})
\end{array}
\end{equation}
\item for $x_{j}=0$ and $\dx _{j}=1$, 
\begin{equation}
\begin{array}{r@{}l}
|x_{j+1}-x_{ex}((j+1)h)| &= O(h^{3})
\\
|\dx_{j+1}-\dx_{ex}((j+1)h)| &= O(ch^{3})
\\
|\mathrm{ME}_{j+1}-\mathrm{ME}_{ex}((j+1)h)| &= O(ch^{3})
\\
|\mathrm{ME}^{c}_{j+1}-\mathrm{ME}^{c}_{ex}((j+1)h)| &= O(ch^{3})
\end{array}
\end{equation}
\end{enumerate}
\item For $x_{j}=\dx_{j}=0$ and $f$ piece-wise exponential of the form $f=f_{j}e^{-c\.}$
in each interval $[jh,(j+1)h]$, with $f_{j}\in \R$,:
\begin{equation}
\begin{array}{r@{}l}
|x_{j+1}-x_{ex}((j+1)h)| &= O(ch^{3})
\\
|\dx_{j+1}-\dx_{ex}((j+1)h)| &= O(ch^{3})
\\
|\mathrm{ME}_{j+1}-\mathrm{ME}_{ex}((j+1)h)| &= O(h^{4})
\\
|\mathrm{ME}^{c}_{j+1}-\mathrm{ME}^{c}_{ex}((j+1)h)| &= O(h^{4})
\end{array}
\end{equation}
The constant in the big-O notation depends on $\max |f_{j+1}|$.
\end{enumerate}
where $x_{ex}$ represents the exact solution of the corresponding SDOF problem.
\end{theorem}
In the notation of the theorem,
\begin{equation}
\mathrm{ME}^{c}_{j} = \frac{1}{2}e^{cjh}kx_{j}^{2} + \frac{1}{2}e^{cjh}\dx_{j}^{2}
\end{equation}
is the modified mechanical energy of the approximate solution, and
$\mathrm{ME}^{c}_{ex}(jh)$ is the modified mechanical energy of the exact solution at
time $t=jh$. The modification is by the exponential factor appearing in the definition
of the relevant quantities, see e.g. eq. \eqref{eqmatrix B}.

\begin{proof}
The iteration is just a restatement of the calculations preceding the statement of the
theorem.

Regarding the error rates, they follow from the calculation of the Taylor expansions 
of the iteration formulas for the corresponding cases. Throughout the proof, $k$ has
been set to $1$ for definiteness. The general case follows by rescaling $t$ accordingly.
The notation shortcut $\Xi _{i} = \Xi_{i}(t,c,1)$ will be used, while it is reminded that
$\w _{d} = \sqrt{1-(c/2)^{2}}$.

Firstly, consider the homogeneous case, $f \equiv 0$.

More precisely, consider the function
\begin{equation}
-\frac{\Xi_{1}+\Xi_{2} }{\Xi_{3}} - e^{-ct/2}(\cos(\w _{d}t)+
\frac{c}{2\w _{d}}\sin(\w _{d} t))
\end{equation}
which is equal to the approximate minus the exact solution for $x_{0}=1$ and
$\dx_{0}=0$. Its Taylor expansion around $t=0$ is
\begin{equation} 
\frac{ct^{3}}{12} + \frac{(4-7c^{2})t^{4}}{120} + O(t^{5})
\end{equation}
Similarly, consider the function
\begin{equation}
-\frac{t}{2}\frac{\Xi_{2} }{\Xi_{3}} -
e^{-ct/2}
\frac{1}{\w _{d}}\sin(\w _{d} t)
\end{equation}
which is the approximate displacement minus the exact one for $x_{0}=0$ and
$\dx_{0}=1$. The expansion of this function reads
\begin{equation}
-\frac{(1-c^{2})t^{3}}{12} + \frac{c(11-7c^{2})t^{4}}{120} + O(t^{5})
\end{equation}
These Taylor expansions account for the good approximation of the displacement in
a homogeneous problem, with $0$ excitation force and non-trivial initial conditions
(displacement and velocity, respectively).

Similarly, consider the function
\begin{equation}
-\frac{2}{t}\frac{\Xi_{1}+\Xi_{2}+\Xi_{3}}{\Xi_{3}} +
e^{-ct/2}
\frac{1}{\w _{d}}\sin(\w _{d} t)
\end{equation}
which is the approximate velocity minus the exact one for $x_{0}=1$ and $\dx_{0}=0$.
The expansion of this function reads
\begin{equation}
-\frac{(1+2c^{2})t^{3}}{60} - \frac{c(5-4c^{2})t^{4}}{120} + O(t^{5})
\end{equation}
Finally, the function
\begin{equation}
-\frac{\Xi_{2}+\Xi_{3}}{\Xi_{3}}
- e^{-ct/2}(\cos(\w _{d}t)-
\frac{c}{2\w _{d}}\sin(\w _{d} t))
\end{equation}
is the approximate velocity minus the exact one for $x_{0}=0$ and $\dx_{0}=1$.
The expansion of this function reads
\begin{equation}
\frac{(1-2c^{2})t^{3}}{60} + \frac{(4-9c^{2}+4c^{4})t^{4}}{120} + O(ct^{5})
\end{equation}

Regarding the conservation of Mechanical Energy for $x_{0}=1$ and $\dx_{0}=0$,
the Taylor expansion of the relevant function (which is too cumbersome to write down)
is
\begin{equation}
\frac{ct^{3}}{12} + \frac{(2-c^{2})t^{4}}{40} + O(ct^{5})
\end{equation}
while that of the modified Mechanical Energy reads
\begin{equation}
\frac{ct^{3}}{12} + \frac{(6+7c^{2})t^{4}}{120} + O(ct^{5})
\end{equation}

Turning to the conservation of Mechanical Energy for $x_{0}=0$ and $\dx_{0}=1$,
the Taylor expansion of the relevant function (again, too cumbersome to write down)
is
\begin{equation}
\frac{c(1-2c^{2})t^{3}}{60} - \frac{(6+c^{2}-8c^{4})t^{4}}{40} + O(ct^{5})
\end{equation}
while that of the modified Mechanical Energy reads
\begin{equation}
\frac{c(1-2c^{2})t^{3}}{60} - \frac{(6-c^{2}-4c^{4})t^{4}}{40} + O(ct^{5})
\end{equation}
These approximations establish the good energy-conservation properties of the
simplest case of the algorithm proposed here-in in the homogeneous case.

Concerning the terms containing the external force, the following can be obtained
by direct calculation. One gets for $F(t) = \int _{0}^{t} e^{ct} f(t)  \dd \t$,
that $F$, as in the undamped case, is projected onto linear polynomials, resulting in
$f$ being of the form $Ke^{-c\.}$. Equivalently,
$\langle f(\. ) , B_{2,3,h}(\. ) \rangle _{c}$ is the projection of $f \in H^{-1}$ onto
functions of the same form $Ke^{-ct}$. Replacing $f(\.) = e^{-c\.}$ in
$\langle f(\. ) , B_{2,3,h}(\. ) \rangle _{c}$ gives
\begin{equation}
\langle
e^{-c\.} , B_{2,3,h}(\. )
\rangle _{c}
=
\int _{0}^{h}
B_{2,3,h}(t) \dd t
= \frac{h}{3}
\end{equation}
i.e. the same result as in the homogeneous case.

The solution of the undamped SDOF system for $k=1$, homogeneous initial conditions and
$f(\.) = e^{-c\.}$ reads
\begin{equation}
x(t) = e^{-ct} - e^{-ct/2}(\cos(\w _{d}t)-
\frac{c}{2\w _{d}}\sin(\w _{d} t))
\end{equation}
The approximate solution in this case yields
\begin{equation}
u_{3} = \frac{h}{3\Xi_{3}}
\end{equation}
The Taylor development around $0$ of the exact minus the approximate value for
the displacement reads
\begin{equation}
-\frac{ct^{3}}{6} - \frac{(2-9c^{2})t^{4}}{60} + O(ct^{5})
\end{equation}
while for the velocity it reads
\begin{equation}
\frac{(1+3c^{2})t^{4}}{60} + \frac{c(9-8c^{2})t^{4}}{120} + O(t^{5})
\end{equation}

Finally, the Taylor development of the difference in the Mechanical Energy reads
\begin{equation}
\frac{(1+3c^{2})t^{4}}{60} + \frac{c(3+14c^{2})t^{5}}{120} + O(t^{6})
\end{equation}
while that of the modified Mechanical Energy reads
\begin{equation}
\frac{(1+3c^{2})t^{4}}{60} + \frac{c(1+8c^{2})t^{5}}{120} + O(t^{6})
\end{equation}
\end{proof}

\section{Error estimates} \label{secerror estimates}

The first step in establishing the estimates for the error of approximation is to obtain
relevant density results for the spaces onto which the excitation force is projected.

\subsection{Some approximation properties}

Firstly, the spaces relevant to approximation, cf. eq.(103)
of \cite{NKSDOFI},
need to be defined in the present context:
\begin{equation} \label{eqdef force spaces bern pol}
\begin{array}{r@{}l}
\FFF = \FFF(p,h) &= \FF ( \PP ) \\
\FFF _{00} = \FFF _{00} (p,h) &= \FF ( \PP _{00}) \\
\FFF_{0}^{0} =\FFF\h (p,h) &= \FF	 ( \PP \h ) \\
\FFF _{ic} &= \FF	 ( \PP _{ic} )
\end{array}
\end{equation}
where
\begin{equation}
\begin{array}{r@{}l}
\PP =  \PP (h) &= \mathrm{vec}(\{ B_{i} \}_{i=1}^{p})\\
\PP _{00}= \PP _{00} (h)&= \mathrm{vec}(\{ B_{i} \}_{i=3}^{p}) \\
\PP \h = \PP \h (h)&= \mathrm{vec}(\{ B_{i} \}_{i=2}^{p-1})\\
\dot{\PP}\h= \dot{\PP}\h (p,h)&= \{ \dot{\f}, \f \in
\PP\h(p,h) 
\}
\\
&= \{
\phi \in \PP(p-1,h),\int \phi = 0
\}\\
\PP _{ic}= \PP _{ic} (h)&= \mathrm{vec}(\{ B_{i} \}_{i=1}^{2})
\end{array}
\end{equation}

The algorithm proposed herein can be summed up as follows. Let a time-step $h$,
a polynomial degree of approximation $p-1 \in \N$, $p\geq 3$, and an interval $[0,\T]$
with $\T = lh$ with $l \in \N ^{*}$ be given. Then, the approximate solution of
\S 9 of \cite{NKSDOFI} is constructed in each subinterval $[(j-1)h,jh]$,
$1\leq j \leq l$ of length $h$, and the final displacement and velocity at
each step $j$ will be used as initial displacement and velocity for the next step.
Therefore, in order to obtain convergence one can choose between letting
$h \ra 0$, $p \ra \infty$, or eventually both. As a consequence, two types of
density results are needed.


From the Stone–Weierstrass approximation theorem (see \cite{RudMathAn}), the following
follows immediately.
\begin{theorem}
As $p \ra \infty$
\begin{equation}
\begin{array}{r@{}l}
\PP(p,h) &\ra \HH ([0,h]) \\
\dot{\PP}\h(p,h) &\ra \{ \f \in
L^{2}_{0} ([0,h]), \int \f = 0
\}
\end{array}
\end{equation}
in the sense of pointwise convergence.
\end{theorem}
\begin{proof}
The first assertion follows directly from the Stone-Weierstrass theorem, while the second
from the fact that $\dot{\PP}\h(p,h)$ is the space of polynomials of degree $p-2$ with
zero mean value, since they are derivatives of polynomials of degree $p-1$ satisfying
homogeneous boundary conditions.
\end{proof}
\begin{corollary}
For any fixed $f \in H^{-1} ([0,h])$,
\begin{equation}
\|f - \pi_{\PP \h}f\|_{H^{-1}} \ra 0 \text{ as } p \ra \infty
\end{equation}
\end{corollary}

\subsection{The error rate due to the projection error term}

The term $\|f - \pi_{\PP \h}f\|_{H^{-1}}$ appears in corollary 9.6 of
\cite{NKSDOFI}, and controls the convergence of the approximation.
The following two propositions will allow the estimation of the approximation error. Under
the weakest possible assumption that $f \in H^{-1}$, which is the case if $f$ has
Dirac-$\delta$ type forces, no rate of convergence can be obtained. Under the
relatively week assumption that $F$, an integral of $f$, be $H^{1}$, which
covers the case where $f$ Heaviside excitation function, the following results can
be proved. Naturally, additional regularity of $f$ improves the rate of convergence
of the method.

\begin{prop} \label{propapprox Lip}
Let $e^{c\.} f \in H^{s} \ra \R $ be the derivative of $F$, a $H^{s+1}$ function such that
$\int _{0}^{h} F = 0$, where $s \geq -1$. Then,
\begin{equation}
\| e^{c\.} f- \pi _{\PP \h} e^{c\.}  f\|_{H^{-1}} \leq C (p-1)^{-(s+1)} h^{s+1} \| F\|_{H^{s+1}(0,h)}
\end{equation}
\end{prop}

This proposition provides an estimate for the first term in the estimate of corollary
9.6 of \cite{NKSDOFI}, the error due to projection onto a
finite dimensional subsbpace. This error is bounded by
\begin{equation}
    \| e^{c\.}  f- \pi _{\PP \h} e^{c\.}  f\|_{H^{-1}} = O ((p-1)^{-(s+1)} h^{s+1})
\end{equation}
for a given excitation function $f$ as $p \ra \infty$ or $h \ra 0$. In the special case where $f$ is
$C^{\infty}$ smooth, the rate of convergence is (theoretically) faster than any negative
power of the degree of polynomial approximation. Naturally, when the degree is high enough
(around $25$ in the experiments carried out be the authors), numerical instabilities appear
so this is a purely theoretical result.
\begin{proof}
It follows by application of corollary \ref{corproj error} and use of the fact that
\begin{equation}
\| e^{c\.}  f- \pi _{\PP \h} e^{c\.}  f\|_{H^{-1}(0,h)} = \|F - \pi _{\dot{\PP} \h} F\|_{L^{2}_{0}(0,h)}
= \|F - T_{p-1} F\|_{L^{2}_{0}(0,h)}
\end{equation}
in the notation of corollary \ref{corproj error}, where the factor $h^{s}$ has been
incorporated into the constant $C$, since $h$ is kept constant.
\end{proof}

\subsection{The error rate due to the misalignment of subspaces}

In the following, the error rates due to the angles appearing in the
statement of corollary 9.6 of \cite{NKSDOFI}
are obtained. The case of an undamped
system is examined first, since the expressions are much simpler, and
the arguments more transparent.
%
%
In this section, the $L^{2}$ norms of the quantities $F_{i,p,h}$ will be understood
as the homogeneous $L^{2}$ norms, equal to
\begin{equation}
\begin{array}{r@{}l}
\|F_{i,p,h}\|_{L^{2}}^{2}
&=
\int _{0}^{h} \left(
F_{i,p,h} - \frac{1}{h} \int  _{0}^{h} F_{i,p,h}
\right)^{2}
\\
&=
\int _{0}^{h}
F_{i,p,h}^{2}
-
\frac{1}{h}
\left(
\int  _{0}^{h} F_{i,p,h}
\right)^{2}
\end{array}
\end{equation}
which is the relevant norm, since the $F_{i,p,h}$ are defined modulo an integration
constant.

Since, by time-reparametrization, $k$ can be brought to $k=1$, with $c$ and $h$
transforming according to the rule of eq. (42) of
\cite{NKSDOFI}, $k$ will be suppressed for simplicity and the standing
assumption will be $k=1$.

The graphs presented in the following paragraphs are obtained using the data available
in
\href{https://drive.google.com/file/d/1llAm5VxTDw3llUq8g9YdSSsbe1azH27h/view?usp=sharing}{this json file}. Reference to exact values will be avoided for
readability, but the interested reader can access the data, stored in the form of
a dictionary.

\subsubsection{The undamped system, preliminary calculations}

Firstly, concerning the integral of the excitation function in the undamped
case when the displacement function is a Bernstein polynomial, the following
holds.

\begin{lemma} \label{lemforce undamped}
    Let $c = 0 $ and $B_{i,p,h}$, $1\leq i\leq p$, be given. Then, for
    $F_{i,p,h} = \int \FF_{0,k} (B_{i,p,h}) $,
    \begin{equation}
        \begin{array}{r@{}l}
        \int _{0}^{h} F_{i,p,h} ^{2} &= \|\dot{B}_{i,p,h}\| _{0}^{2} +
        k^{2}\|\int B_{i,p,h}\| _{0}^{2} +
        \frac{2h}{p}k
        B_{i,p,h} (h ) -
        2 k \|B_{i,p,h}\| _{0}^{2}
        \\
        &= \|\dot{B}_{i,p,h}\| _{0}^{2} +
        k^{2}\|\int B_{i,p,h}\| _{0}^{2} +
        \frac{2h}{p}k
        \d_{i,p} -
        2 k \|B_{i,p,h}\| _{0}^{2}
        \\
        \int F_{i,p,h} &= \d _{i,p} - \d _{i,1} +
        k\frac{h^{2}}{p}\frac{p-i+1}{p+1}
        \\
        \pi _{(\dot{\PP} \h)^{\perp}} F_{i,p,h} &= k \pi _{(\dot{\PP} \h)^{\perp}}
        \int B_{i,p,h}
        \end{array}
    \end{equation}
\end{lemma}
\begin{proof}
    The first item follows directly from lemma 4.2 of
    \cite{NKSDOFI} for
    $c=0$ and an integration by parts.

    The second item follows from lemmas \ref{lemBern deriv},
    \ref{lemBern integral} and \ref{lemBern integrals}.
    

    The third item follows from the same lemma, using the fact that
    $\dot{B}_{i,p,h} $ has degree $p-2$ just like functions in
    $\dot{\PP} \h$, and functions in $\dot{\PP} \h$ are considered modulo their
    mean value.
\end{proof}

The previous lemma has the following corollary, whose proof follows immediately from
lem. \ref{lemBern integrals}, \ref{lemBern prod} and prop. \ref{propBernLeg}.
\begin{corollary} \label{corforce proj exact}
    It holds true that
    \begin{equation}
    \begin{array}{r@{}l}
        \pi _{(\dot{\PP} \h)^{\perp}} F_{i,p,h} (\. )
        &= \pi _{(\dot{\PP} \h)^{\perp}} k\frac{h}{p}
    \sum\limits _{j=i+1}^{p+1}
    B_{j,p+1,h}(\. )
    \\
    &= k\frac{h}{p} \sum\limits _{m = p-1}^{p} \sqrt{2m+1} 
    \sum\limits _{j=i+1}^{p+1} \Lambda ^{-1} _{jm} (p+1) \tilde{L}_{m}(\. )
    \\
    &= k\frac{h}{p} \sum\limits _{m = p-1}^{p} \sqrt{2m+1} 
    \tilde{F} _{imp} \tilde{L}_{m}(\. )
    \end{array}
    \end{equation}
Consequently,
\begin{equation}
\|\pi _{(\dot{\PP} \h)^{\perp}} F_{i,p,h}\|_{L^{2}}^{2} =
k^{2}\frac{h^{2}}{p^{2}} \sum\limits _{m = p-1}^{p} (2m+1)
    \tilde{F} _{imp}^{2}
\end{equation}
\end{corollary}

The (symmetric, positive-semi-definite) matrix of scalar products
\begin{equation}
\langle
\pi _{\dot{\PP} \h} F_{i,p,h} ,\pi _{\dot{\PP} \h} F_{j,p,h}
\rangle _{0}
\end{equation}
is constructed, and its eigenvectors $\{ \mathbf{e}_j \}_{j=1}^{p}$ calculated.
This is the matrix of scalar products of forces appearing in the construction of
the approximate solution.

The first two eigenvectors correspond to $0$, or very small, eigenvalues. Linear
combinations of these eigenvectors construct the solutions to the homogeneous problem,
with $0$ excitation function and non-trivial initial conditions. The first two
co-ordinates of the null-eigenvectors are linearly independent, allowing for a solution
of the problems
\begin{equation}
\mathbf{A}. \mathbf{\tilde{e}}_1 = 
\begin{bmatrix}
1 \\ 0
\end{bmatrix}
 \text{ and } \mathbf{A}. \mathbf{\tilde{e}}_2 = 
\begin{bmatrix}
0 \\ 1
\end{bmatrix}
\end{equation}
where
\begin{equation}
\mathbf{A} =
\begin{bmatrix}
\mathbf{e}_{1,1} & \mathbf{e}_{2,1}
\\
\mathbf{e}_{1,2} & \mathbf{e}_{2,2}
\end{bmatrix}
\end{equation}

Then, following the corresponding part of \eqref{eq lin sys simple aug}, the eigenvectors
$\mathbf{\tilde{e}}_j, j=1,2$ are expressed in terms of linear combinations of initial
conditions resulting in a basis of the nullspace $\mathbf{e}_{disp}$ and
$\mathbf{e}_{vel}$ corresponding to the solutions with unit initial displacement and
$0$ velocity; and $0$ initial displacement and unit velocity, respectively.

Following the numerical study carried out following prop. \ref{propBernLeg},
the equivalent study of the projection of
\begin{equation}
F \in \mathrm{Vec}(\{ F_{i,p,h}\}_{i=1}^{i=p})
\end{equation}
as in the corollary here above can now be presented.

\subsubsection{The undamped system, $h=T$}

The study is carried out for $h=T$ in a time-scale at which $k=1.$, i.e.
for $h=T = 2\pi $ where $T$ is the natural period of the system. The
scaling properties with respect to $h$ will be studied afterwards.

It should be noted that in the undamped case the terms coming from the derivatives are
in the kernel of the projection operator $\pi _{(\dot{\PP} \h)^{\perp}}$, while the
integrals are not. Their projection on the orthogonal of the kernel, however, remains
very small as the study will establish.

\begin{figure}[h!]
  \centering
  \includegraphics[width=0.7\linewidth]{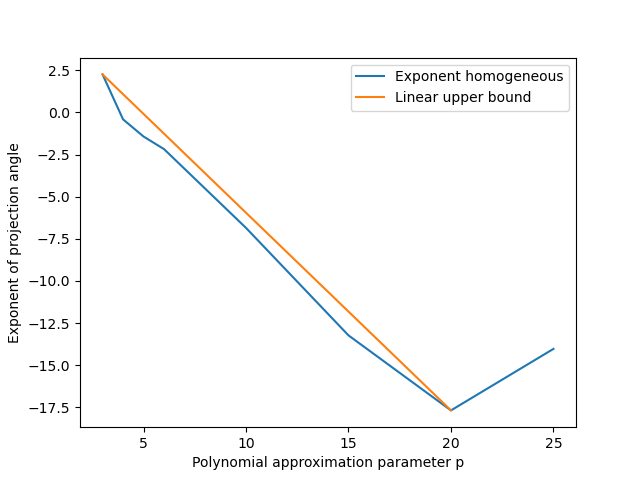}
  \caption{The base-$p$ exponent expressing the error of projection of the force
  corresponding the homogeneous problem. This is the quantity $s_{h}(p)$ in the
  established notation.}
  \label{fig:exp_force_hom}
\end{figure}

The quantity
\begin{equation} \label{eqdef sh}
    s_{h} = s_{h}(p) = \max _{ \max \{ |x_{0}|, |\dx_{0}|\} = 1}
    \frac{\log  \|\pi _{(\dot{\PP} \h)^{\perp}}
    (x_{0}\mathbf{e}_{disp} + \dx_{0}\mathbf{e}_{vel})\|_{L^{2}}}{\log p}
\end{equation}
can be calculated, the subscript $h$ standing for "homogeneous", so that
\begin{equation}
	\|F_{er,h}\|_{L^{2}} =
    \|\pi _{(\dot{\PP} \h)^{\perp}} F_{h}\|_{L^{2}} \leq
    \max \{ |x_{0}|, |\dx_{0}|\} p^{s_{h}(p)}, \forall F \in \mathrm{Vec}(\{ \mathbf{e}_{1}
    ,
    \mathbf{e}_{2}\})
\end{equation}
The full signature of the function actually reads $s_{h}(p,h,c,k=1.)$
so that $s_{h}(p) = s_{h}(p,h=T,c=0.)$ for brevity in notation.
In other terms, the quantity $K \sin \theta _{h}$ of prop. 9.5 of
\cite{NKSDOFI} satisfies
\begin{equation}
K \sin \theta _{h} = p^{s_{h}(p)}
\end{equation}
The graph in fig. \ref{fig:exp_force_hom} plots an estimate of $s_{h}(p)$. The
estimate was obtained by randomly sampling $10000$ initial conditions and keeping
the maximum angle over these samples. All null eigenvalues were smaller than
$\num{1e-2}$ (in fact much smaller, but a threshold was set in order to assure
smallness).
It can be seen that the point of diminishing returns is reached at around $p=20$, where
$K \sin \theta _{h} \approx \num{1e-23}$. For $p=25$ the corresponding value is
$\num{2.5 e-20}$. The rise in the exponent is arguably due to numerical instabilities
and the diminishing returns point could be pushed further by using more sophisticated
numerical analysis machinery.
The broken-line model $s^{bl}_{h}(p)$ is given by
\begin{equation}
s^{bl}_{h}(p) = 
\begin{cases}
-1.1725005559764048 p + 5.771540799001279, 0\leq p \leq 20,
\\
0.729932981104102 p -32.27712994260885, 20 \leq p \leq 25,
\end{cases}
\end{equation}
and it is quite pessimistic in the first leg. In any case, for all relevant values of
$p$, $s_{h}(p) \leq s^{bl}_{h}(p)$.

\begin{figure}[h!]
  \centering
  \includegraphics[width=0.7\linewidth]{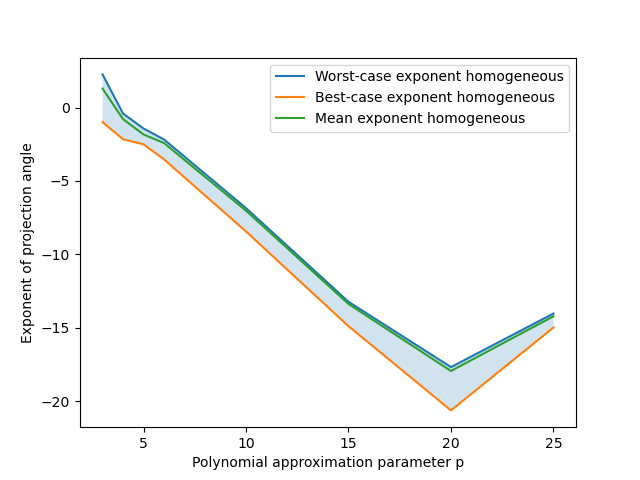}
  \caption{The range of base-$p$ exponents expressing the error of projection of the
  force corresponding the homogeneous problem. The quantities $s_{h}(p)$,
  $\bar{s}_{h}(p)$, and $\underline{s}_{h}(p)$ in the established notation are
  plotted, representing the exponent for the worst-case angle, the mean angle and
  the best-case angle.}
  \label{fig:exp_force_hom_range}
\end{figure}

This is the worst-case scenario for the convergence of the algorithm. A range of observed
rates of convergence can also be visualized by plotting the exponent of the mean value of
the projection angle over the samples, as well as the best-case scenario. The exponent
corresponding to the mean value of the angle is denoted by $\bar{s}_{h}(p)$, and the
best-case scenario by $\underline{s}_{h}(p)$. They are plotted in fig.
\ref{fig:exp_force_hom_range}, where it can be seen that the mean value closely tracks
the worst-case exponent.

A significant spectral gap separates the two null eigenvectors with
the rest of the eigenvectors, whose linear combinations construct the solutions to the
non-homogeneous problem, with non-trivial excitation function. Subtraction of the
corresponding multiples of $\mathbf{\tilde{e}}_j, j=1,2$ from
$\{ \mathbf{e}_j \}_{j=3}^{p}$ imposes $0$ initial conditions, and the eigenvectors
obtained in this way are denoted by $\{ \mathbf{\tilde{e}}_j \}_{j=3}^{p}$. It should
be noted that $\{ \mathbf{\tilde{e}}_j \}_{j=1}^{p}$ is still a basis of eignvectors,
only not orthogonal. Then, the quantity
\begin{equation} \label{eqdef snh}
    s_{nh} = s_{nh}(p) = \max _{\substack{ F \in \mathrm{Vec} \{ \mathbf{\tilde{e}}_j \}_{j=3}^{p}
    \\
    \|F \|_{L^{2}} \leq 1}}
    \frac{\log  \|\pi _{(\dot{\PP} \h)^{\perp}} F\|_{L^{2}}}{\log p}
\end{equation}
can be calculated, so that
\begin{equation}
    \|\pi _{(\dot{\PP} \h)^{\perp}} F\|_{L^{2}} \leq
    \| F\|_{L^{2}} p^{s_{nh}(p)}, \forall F \in \mathrm{Vec}\{ \mathbf{\tilde{e}}_j \}_{j=3}^{p}
\end{equation}
so that the angle $\theta _{nh}$ of prop. 9.5 of
\cite{NKSDOFI} satisfies
\begin{equation}
\sin \theta _{nh} = p^{s_{h}(p)}
\end{equation}
The sine of the angle, instead of the tangent, is calculated for purely practical
reasons, as dividing with the hypotenuse is numerically more stable. In the limit
$\theta _{nh} \ra 0$, which is the interesting case, the quantities are equal for
all practical reasons. For a given approximate solution $x_{ap}$ in
$\mathrm{Vec}\{ B_{i,p,h} \}_{i=3}^{p}$, the quantity $p^{s_{h}(p)}$ expresses the
proportion of the force $\FF (x_{ap})$ that contributes to the error term.

\begin{figure}[h!]
  \centering
  \includegraphics[width=0.7\linewidth]{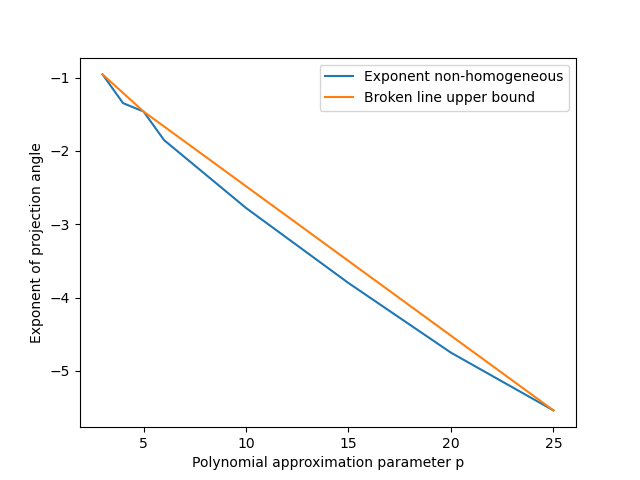}
  \caption{The base-$p$ exponent expressing the error of projection of the force
  corresponding the homogeneous problem. This is the quantity $s_{h}(p)$ in the
  established notation.}
  \label{fig:exp_force_non_hom}
\end{figure}

The graph in fig. \ref{fig:exp_force_non_hom} plots an estimate of $s_{nh}(p)$. The
estimate was obtained by randomly sampling $10000$ linear combinations of vectors
in $\{ \mathbf{e}_j \}_{j=3}^{p}$ and getting the maximum angle over these samples.
The spectral gap condition was that the minimal non-null eigenvalue be at least
$\num{1e2}$ times greater than the biggest null eigenvalue. The spectral gap is actually
a lot bigger than $2$ orders of magnitude, but, as before, a threshold was set in order
to assure a sufficient spectral gap.
In the non-homogeneous problem, no point of diminishing returns is observed, which
corroborates the hypothesis of not sufficiently strong numerical analysis machinery
for treating the null eigenvalues that are (very close to) $0$.
The linear model $s^{\ell}_{nh}(p)$ is given by
\begin{equation}
s^{bl}_{nh}(p) = 
\begin{cases}
 -0.2530595464212411 p -0.19694374410552673, 0\leq p \leq 5,
\\
-0.20385935288152152 p -0.44294471180412476, 5 \leq p \leq 25,
\end{cases}
\end{equation}
and it is also a bit pessimistic in the first leg. In any case, for all relevant values
of $p$, one still has that $s_{nh}(p) \leq s^{bl}_{nh}(p)$. It can be remarked that,
already for $p=3$ the exponent is negative, unlike for the homogeneous problem. However,
the factor is merely of the order of $1/3$, not enough to produce an acceptable error
rate.

\begin{figure}[h!]
  \centering
  \includegraphics[width=0.7\linewidth]{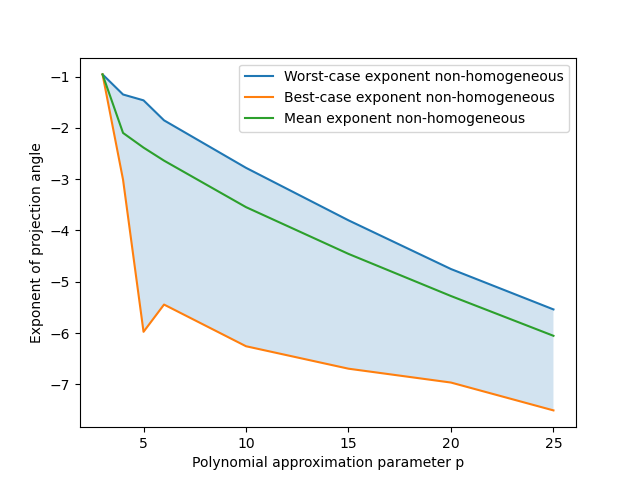}
  \caption{The range of base-$p$ exponents expressing the error of projection of the
  force corresponding the homogeneous problem. The quantities $s_{nh}(p)$,
  $\bar{s}_{nh}(p)$, and $\underline{s}_{nh}(p)$ in the established notation are
  plotted, representing the exponent for the worst-case angle, the mean angle and
  the best-case angle.}
  \label{fig:exp_force_non_hom_range}
\end{figure}

As for the homogeneous problem, the best-case scenario angle and the exponent of the mean
angle are plotted in fig. \ref{fig:exp_force_non_hom_range}, showing that the worst-case
scenario tends to be more pessimistic in the homogeneous than in the non-homogeneous
case, with the exponent of the mean angle sitting farther than the worst case angle.

\subsubsection{The undamped system, dependence on $h$}

In this section, the dependence of the angles of projection on the timestep parameter
$h$ will be studied. The scaling of all quantities involved is highly non-trivial,
since the derivative scales like $h^{-1}$, the integral like $h$, and the products
of derivatives and integrals are constant with respect to $h$.

It should be expected that for small $h$ the dominant term come from the derivative,
leading to small angles since $\pi _{(\dot{\PP} \h)^{\perp}} \dot{B}_{i,p,h} \equiv 0$.
For intermediate values of $h$, all terms should contribute roughly in par, while for
larger $h$ the integral should start to dominate, leading to greater angles, since
integrals are not in the kernel of the projection operator
$\pi _{(\dot{\PP} \h)^{\perp}}$.
The products of derivatives and integrals have complicated contribution to the norm
of the force, since the derivative can have either sign, while the integrals are
always $\geq 0$.

Following eq. \eqref{eqdef sh}, the quantity
\begin{equation} \label{eqdef sh_h}
    s_{h} = s_{h}(p,h) = \max _{ \max \{ |x_{0}|, |\dx_{0}|\}= 1}
    \frac{\log  \|\pi _{(\dot{\PP} \h)^{\perp}}
    (x_{0}\mathbf{e}_{disp} + \dx_{0}\mathbf{e}_{vel})\|_{L^{2}}
     - s_{h}(p,T)\log p}{\log (h/T)}
\end{equation}
can be calculated, so that
\begin{equation} \label{eqconv rage hom ph}
    \|\pi _{(\dot{\PP} \h)^{\perp}} F\|_{L^{2}} \leq
    ( \max |x_{0}|, |\dx_{0}|) p^{s_{h}(p,T)}
    \left(
    \frac{h}{T}
    \right)^{s_{h}(p,h)}, \forall F \in \mathrm{Vec}(\{ \mathbf{e}_{1}
    ,
    \mathbf{e}_{2}\})
\end{equation}

As for $h=T$, the same number of points are sampled, keeping the worst exponent
The same thresholds for null eigenvalues an the spectral gap. Timesteps were sampled
in the range from $0.001T$ up to $4T$, i.e. in the range of $ 0.1 \%$ to $400\%$ of 
the natural period of the system with $k=1$ and $c=0$. 

\begin{figure}[h!]
  \centering
  \includegraphics[width=0.7\linewidth]{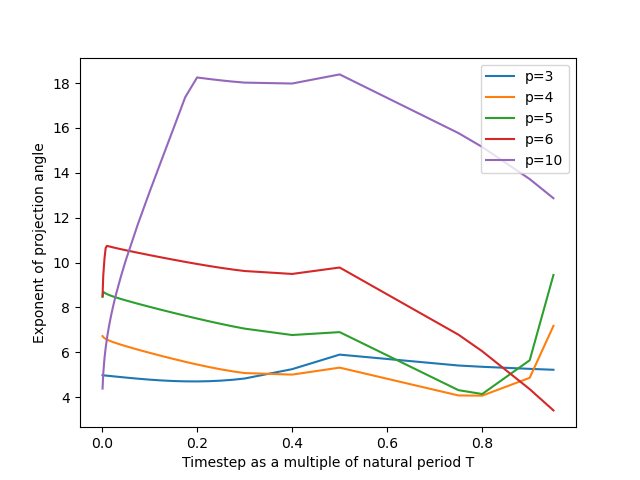}
   \vspace{-6mm}
  \caption{The base-$h$ exponent expressing the error of projection of the force
  corresponding the homogeneous problem, for timestep smaller than the natural period
  and $3 \leq p \leq 10$.
  This is the quantity $s_{h}(p,h)$ in the established notation.}
  \label{fig:exp_force_h_hom_low_small}
\end{figure}

\begin{figure}[h!]
  \centering
  \includegraphics[width=0.7\linewidth]{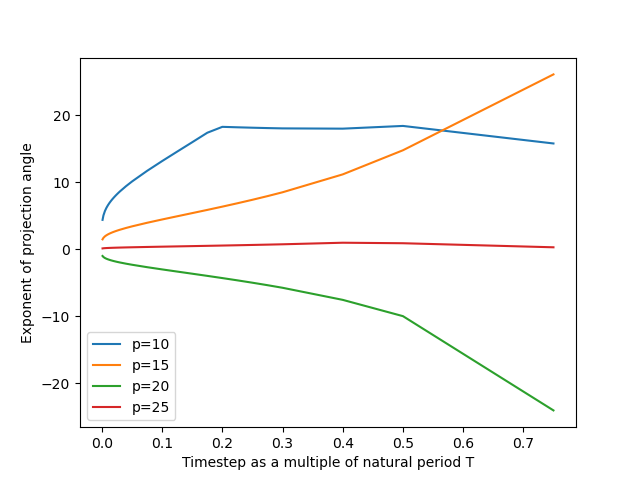}
   \vspace{-6mm}
  \caption{The base-$h$ exponent expressing the error of projection of the force
  corresponding the homogeneous problem, for timesteps smaller than the natural period,
  and $10 \leq p \leq 25$. This is the quantity $s_{h}(p,h)$ in the established
  notation.}
  \label{fig:exp_force_h_hom_low_large}
\end{figure}

\begin{figure}[h!]
  \centering
  \includegraphics[width=0.7\linewidth]{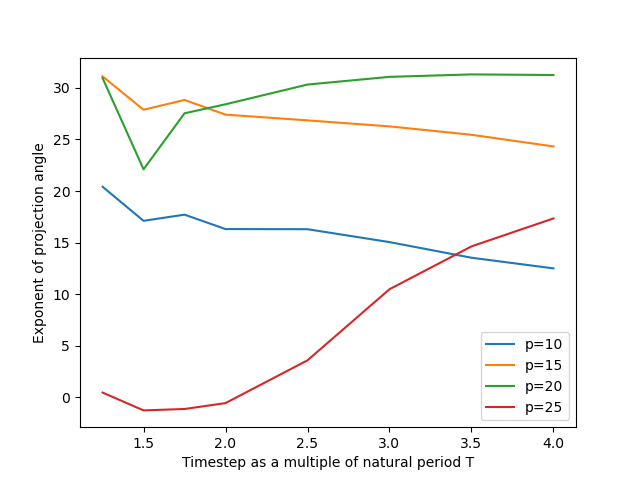}
   \vspace{-6mm}
  \caption{The base-$h$ exponent expressing the error of projection of the force
  corresponding the homogeneous problem, for timesteps greater than the natural period
  and $10 \leq p \leq 25$.
  This is the quantity $s_{h}(p,h)$ in the established notation.}
  \label{fig:exp_force_h_hom_high}
\end{figure}

The function $s_{h}(p,h)$ is plotted in figs. \ref{fig:exp_force_h_hom_low_small}
and  \ref{fig:exp_force_h_hom_low_large} for values of $h<T$ and small, resp. large
values of $p$, and in  fig. \ref{fig:exp_force_h_hom_high} for values of $h>T$ and
large values of $p$. The split in $h$ was necessary in order to avoid division by $0$
at $h=T$, and the split in $p$ for clarity.

The graph of fig. \ref{fig:exp_force_h_hom_low_small}, featuring values of $p$ in the
range $3\leq p \leq 10$, shows a very fast rate of
convergence, as all exponents are consistently $>4$, which is an added factor of
convergence on top of the factor $p^{s_{h}(p,T)}$, c.f. eq. \eqref{eqconv rage hom ph}.
The graph of fig. \ref{fig:exp_force_h_hom_low_large}, featuring $10\leq p\leq 25$,
shows a more complicated behavior. While for $p=10,15$ the exponent is positive,
leading to a
faster convergence, smaller $h$ slows down convergence for $p=20$. The latter remains
fast, nonetheless, since, say
\begin{equation}
K \sin \theta _{h} \leq 0.5^{-10} 20 ^{s^{bl}_{h}(20)} \approx \num{1e-20}.
\end{equation}
The convergence for $p=25$ is rather indifferent to $h$, if only slightly accelerated
by small values of $h$, with exponents in the range $(0,1)$.

The graph of fig. \ref{fig:exp_force_h_hom_high} shows equally a non-monotonous behavior
of the exponent $s_{h}(p,h)$ for $h>T$. The graph features only high values of $p$,
namely $10 \leq p \leq 25$, since for low values of $p$ timesteps greater than the
natural period $T$ of the system will not produce a small error rate. Already for $p=15$
convergence can remain fast for large values of $h$, since for $h=2T$,
\begin{equation}
K \sin \theta _{h} \leq 2^{-30} 15 ^{s^{bl}_{h}(15)} \approx \num{1e-5}.
\end{equation}
Convergence is faster for $p=20$, while the dependence on $h$ gets milder for $p=25$.

Concerning the non-homogeneous problem, following eq. \eqref{eqdef snh}, the quantity
\begin{equation} \label{eqdef sh_nh}
    s_{nh} = s_{nh}(p,h) = \max _{\substack{ F \in \mathrm{Vec} \{ \mathbf{\tilde{e}}_j \}_{j=3}^{p}
    \\
    \|F \|_{L^{2}} \leq 1}}
    \frac{\log  \|\pi _{(\dot{\PP} \h)^{\perp}} F\|_{L^{2}} - s_{nh}(p,1) \log p}{\log h}
\end{equation}
can be calculated, so that
\begin{equation}
    \|\pi _{(\dot{\PP} \h)^{\perp}} F\|_{L^{2}} \leq
    \| F\|_{L^{2}} p^{s_{nh}(p,1)} h^{s_{nh}(p,h)}, \forall F \in \mathrm{Vec}\{ \mathbf{\tilde{e}}_j \}_{j=3}^{p}
\end{equation}

\begin{figure}[h!]
  \centering
  \includegraphics[width=0.7\linewidth]{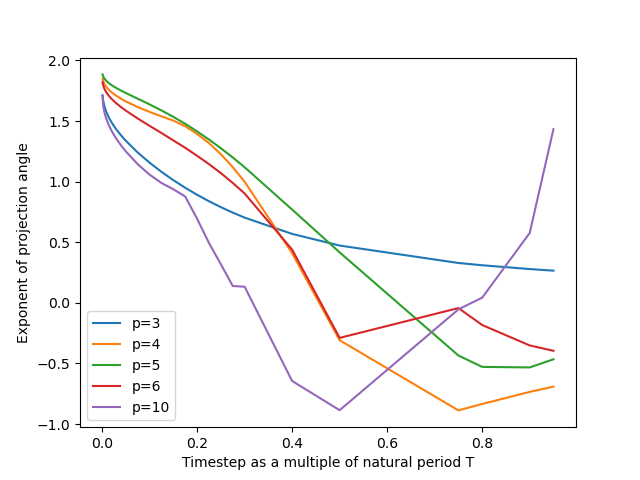}
   \vspace{-6mm}
  \caption{The base-$h$ exponent expressing the error of projection of the force
  corresponding the non-homogeneous problem, for timesteps smaller than the natural
  period, and $3 \leq p \leq 10$. This is the quantity $s_{nh}(p,h)$ in the established
  notation.}
  \label{fig:exp_force_h_non_hom_small}
\end{figure}

\begin{figure}[h!]
  \centering
  \includegraphics[width=0.7\linewidth]{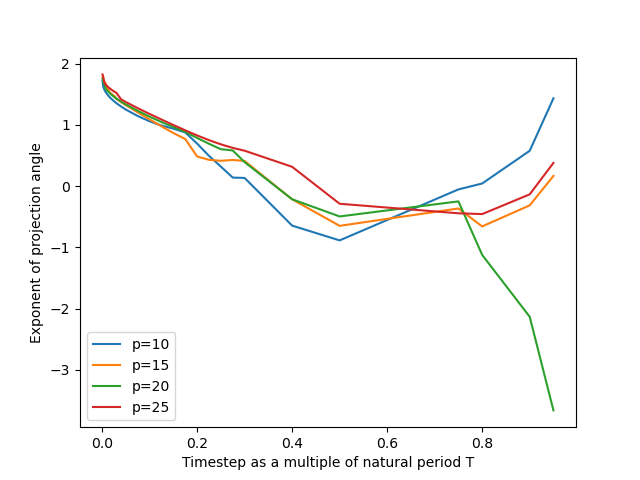}
   \vspace{-6mm}
  \caption{The base-$h$ exponent expressing the error of projection of the force
  corresponding the non-homogeneous problem, for timesteps smaller than the natural
  period, and $10 \leq p \leq 25$. This is the quantity $s_{nh}(p,h)$ in the established
  notation.}
  \label{fig:exp_force_h_non_hom_large}
\end{figure}

\begin{figure}[h!]
  \centering
  \includegraphics[width=0.7\linewidth]{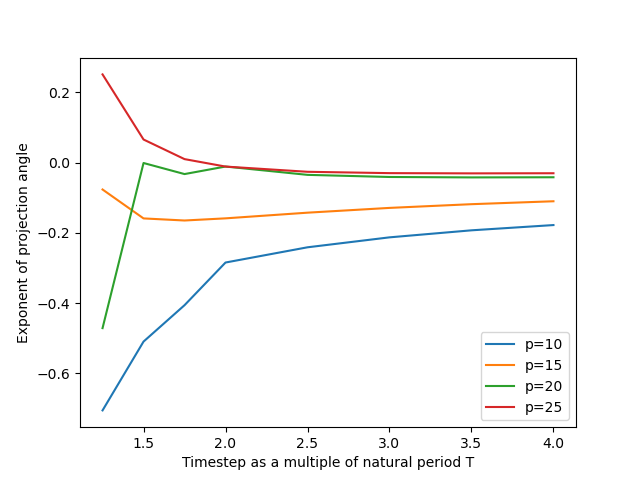}
   \vspace{-6mm}
  \caption{The base-$h$ exponent expressing the error of projection of the force
  corresponding the non-homogeneous problem, for timesteps greater than the natural period
  and $10 \leq p \leq 25$.
  This is the quantity $s_{nh}(p,h)$ in the established notation.}
  \label{fig:exp_force_h_non_hom_high}
\end{figure}

The function $s_{nh}(p,h)$ is plotted in figs. \ref{fig:exp_force_h_non_hom_small},
\ref{fig:exp_force_h_non_hom_large}, and \ref{fig:exp_force_h_non_hom_high}, with the
same parameters as for the homogeneous problem.

As show fig. \ref{fig:exp_force_h_non_hom_small}, for low-order polynomial approximation
and $h$ smaller than $20\%$ of the natural period, the smallness of the timestep
contributes to the rate of convergence with an exponent
ranging between linear and quadratic, and getting better as $h\ra 0$. As $h$ becomes
comparable with the natural period of the system, the corresponding factor never eats
too much on the factor produced by $p^{s_{nh}(p)}$, but ceases to contribute.

The same picture holds for higher order polynomial approximation, if only for slightly
larger values of $h$, as can be seen in fig. \ref{fig:exp_force_h_non_hom_large}. It
is reminded that the factor $p^{s_{nh}(p)}$ is, in this case, orders of magnitute
smaller, so that the effect of the negative exponent for $h/T$ close to $1^{-}$ is
negligible.

Moving on to timesteps greater than the natural period of the system, it can be seen
in fig. \ref{fig:exp_force_h_non_hom_high} that for $p=10,15$ the factor $h/T$
does contribute to the convergence rate, since the exponent is negative. Convergence
becomes practically indifferent to the value of the timestep for $p=20,25$ since the
exponent remains negative, but is very close to $0$.

These graphs suggest that mixed approximation schemes, using different degrees of
polynomial approximation for the homogeneous and non-homogeneous part, are conceivable,
but are outside the scope of the present work.

\subsubsection{The damped system, preliminary calculations}

Concerning the integral of the excitation function in the damped case when the
displacement function is a Bernstein polynomial, the following can be obtained.

\begin{lemma} \label{lemforce damped}
    Let $c \geq 0 $ and $B_{i,p,h}$, $1\leq i\leq p$, be given. Then, for
    $F_{i,p,h} = \int e^{c\.} \FF_{0,k} (B_{i,p,h}(\.)) $,
    \begin{equation}
        \begin{array}{r@{}l}
        \int _{0}^{h} F_{i,p,h} ^{2} &= \|\dot{B}_{i,p,h}\| _{c}^{2} +
        k^{2}\|\int B_{i,p,h}\| _{c}^{2} + 
        2k \langle \dot{B}_{i,p,h} , \int e^{c\.}B_{i,p,h}
        \rangle _{c}
        \\
        \int F_{i,p,h} &= \int _{0}^{h}e^{c\.} \dot{B}_{i,p,h}(\.) +
        k\int _{0}^{h}e^{c\.} B_{i,p,h}(\.)
        \\
        \pi _{(\dot{\PP} \h)^{\perp}} F_{i,p,h} &= 
		\pi _{(\dot{\PP} \h)^{\perp}}e^{c\.} \dot{B}_{i,p,h}(\.)
		+        
        k \pi _{(\dot{\PP} \h)^{\perp}}
        \int e^{c\.} B_{i,p,h}
        \end{array}
    \end{equation}
\end{lemma}
The lemma just states the definitions of the relevant quantities. Closed forms involving
the hypergeometric function can be obtained, using lemmas \ref{lemBern integral c}
and \ref{lemBern times exp}, and cor. \ref{corBern integral c}. They are, however,
very cumbersome and not easy to manipulate.

The statement corresponding to cor. \ref{corforce proj exact} also becomes difficult
to write down. What can be said about the projection is that multiplication by the
exponential leaks mass from the derivative into the orthogonal of the kernel of the
projection operator $\pi _{(\dot{\PP} \h)^{\perp}}$, weighed by powers of $ch$.

At the linear level, $e^{ct} = 1 + ct + O(c^{2}t^{2})$, and this results in
\begin{equation}
\pi _{(\dot{\PP} \h)^{\perp}} c t \dot{B}_{i,p,h}(t) = 
c \pi _{(\dot{\PP} \h)^{\perp}} t \dot{B}_{i,p,h}(t)
\end{equation}
being non-trivial, since $t \dot{B}_{i,p,h}(t)$ is of order $p-1$, c.f. also lem.
\ref{lemBer momentum}.
In the same way, the presence of the exponential factor leaks more $L^{2}$ mass
of the integral $\int e^{c\.} B_{i,p,h}$ from the kernel of the projection operator to
its orthogonal.
In both cases, $L^{2}$ mass is added both in the kernel and its orthogonal, and the
trade-off seems impossible to analyze accurately.

The same procedure as in the undamped case can be followed in order to obtain numerical
results concerning the matrix of the products of forces.

\subsubsection{The effect of damping, $h=T$}

The first part of the study is performed by keeping $h=T$, constant, and studying
how the curve of \ref{fig:exp_force_hom} changes as damping is switched on.
As explained in the previous paragraph, the dependence of the involved quantities
on $c$ are very complicated, so, once again the authors are restricted to a
numerical study of the phenomenon.

\begin{figure}[h!]
  \centering
  \includegraphics[width=0.7\linewidth]{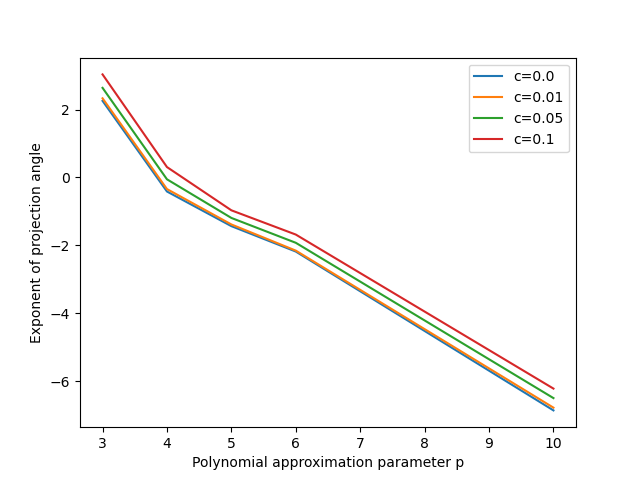}
   \vspace{-6mm}
  \caption{The base-$p$ exponent expressing the error of projection of the force
  corresponding the homogeneous problem. The graph features low values of 
  $p$, $3\leq p \leq 10$, and low values of $c$, $0 \leq c \leq 0.1$. This is the
  quantity $s_{h}(p,c,T)$ in the established notation.}
  \label{fig:exp_force_p_c_hom_low_low}
\end{figure}

\begin{figure}[h!]
  \centering
  \includegraphics[width=0.7\linewidth]{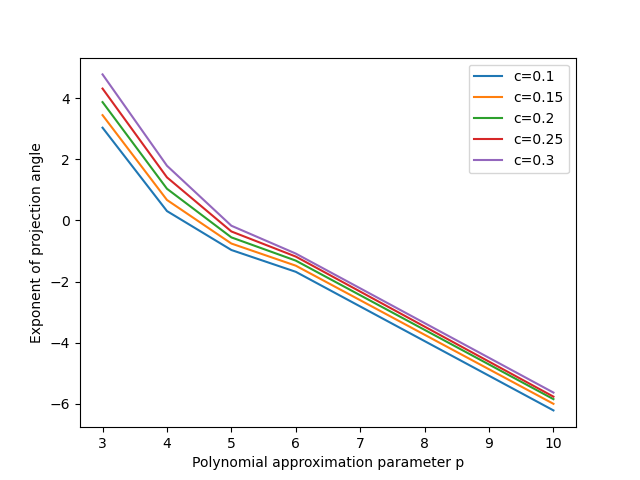}
   \vspace{-6mm}
  \caption{The base-$p$ exponent expressing the error of projection of the force
  corresponding the homogeneous problem. The graph features low values of 
  $p$, $3\leq p \leq 10$, and high values of $c$, $0.1 \leq c \leq 0.3$. This is the
  quantity $s_{h}(p,c,T)$ in the established notation.}
  \label{fig:exp_force_p_c_hom_low_high}
\end{figure}

\begin{figure}[h!]
  \centering
  \includegraphics[width=0.7\linewidth]{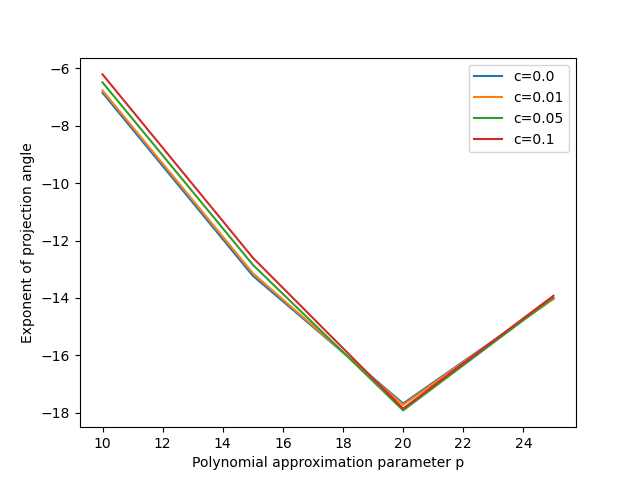}
   \vspace{-6mm}
  \caption{The base-$p$ exponent expressing the error of projection of the force
  corresponding the homogeneous problem. The graph features high values of 
  $p$, $10\leq p \leq 25$, and low values of $c$, $0 \leq c \leq 0.1$. This is the
  quantity $s_{h}(p,c,T)$ in the established notation.}
  \label{fig:exp_force_p_c_hom_high_low}
\end{figure}

\begin{figure}[h!]
  \centering
  \includegraphics[width=0.7\linewidth]{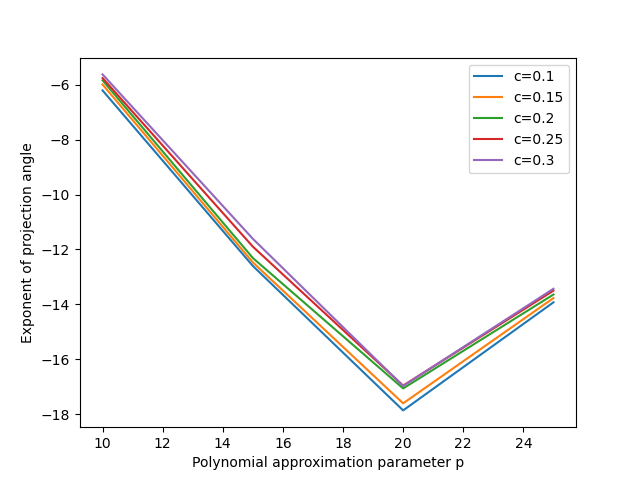}
   \vspace{-6mm}
  \caption{The base-$p$ exponent expressing the error of projection of the force
  corresponding the homogeneous problem. The graph features high values of 
  $p$, $10\leq p \leq 25$, and high values of $c$, $0.1 \leq c \leq 0.3$. This is the
  quantity $s_{h}(p,c,T)$ in the established notation.}
  \label{fig:exp_force_p_c_hom_high_high}
\end{figure}

To this end and regarding the homogeneous problem, the quantity of eq. \eqref{eqdef sh_h}
is calculated, taking into account the dependence of $s_{h} = s_{h}(p,c,T)$ on the
damping coefficient $c$. The result is plotted in figs.
\ref{fig:exp_force_p_c_hom_low_low}, \ref{fig:exp_force_p_c_hom_low_high},
\ref{fig:exp_force_p_c_hom_high_low}, and \ref{fig:exp_force_p_c_hom_high_high}. The
graphs are split for $p$ and for $c$ for clarity.

The figures show that, in general, the projection angle deteriorates as $c$ grows.
As in the undamped case, the part of the graph that is monotonously increasing for $p$
with $c$ held constant can probably be attributed to issues in the numerical calculation
of null eigenvectors.

The dominant effect of damping is to increase the mass in the kernel of the projection
operator faster than in the orthogonal of the kernel. The effect becomes weaker as
$p$ increases, leading to the inverse phenomenon for $p=20$ and small values of $c$
as can be seen in graph \ref{fig:exp_force_p_c_hom_high_low}. 

Overall, the introduction of damping slows down the rate of convergence, since for all
$p$ it holds that
\begin{equation}
s_{h}(p,0.3,T) \leq s_{h}(p,0,T) + 2.55
\end{equation}
with the losses being bigger for small $p$.

\begin{figure}[h!]
  \centering
  \includegraphics[width=0.7\linewidth]{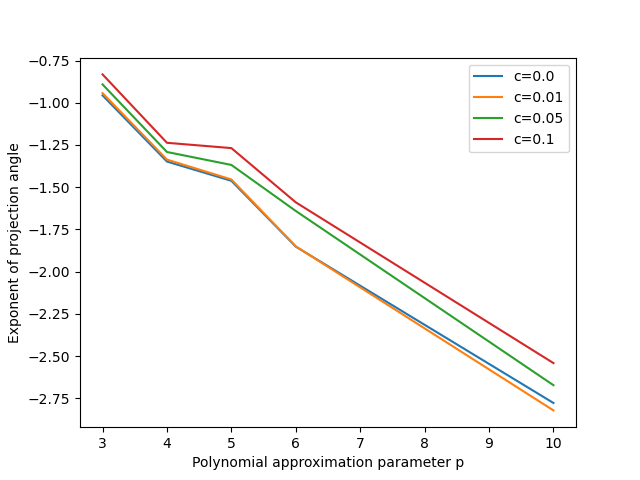}
   \vspace{-6mm}
  \caption{The base-$p$ exponent expressing the error of projection of the force
  corresponding the non-homogeneous problem. The graph features low values of 
  $p$, $3\leq p \leq 10$, and low values of $c$, $0 \leq c \leq 0.1$. This is the
  quantity $s_{nh}(p,c,T)$ in the established notation.}
  \label{fig:exp_force_p_c_non_hom_low_low}
\end{figure}

\begin{figure}[h!]
  \centering
  \includegraphics[width=0.7\linewidth]{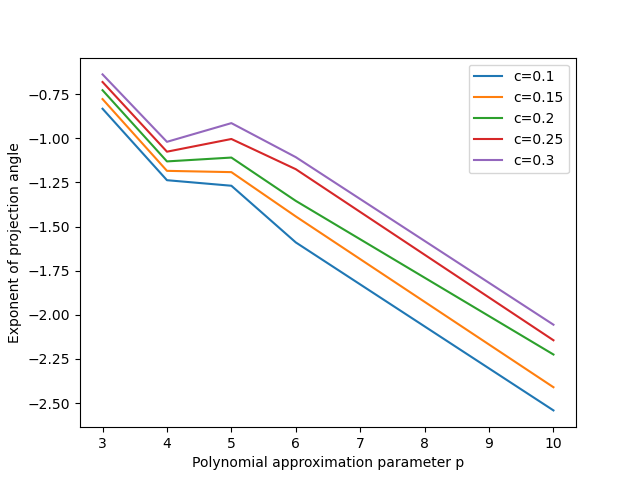}
   \vspace{-6mm}
  \caption{The base-$p$ exponent expressing the error of projection of the force
  corresponding the non-homogeneous problem. The graph features low values of 
  $p$, $3\leq p \leq 10$, and high values of $c$, $0.1 \leq c \leq 0.3$. This is the
  quantity $s_{nh}(p,c,T)$ in the established notation.}
  \label{fig:exp_force_p_c_non_hom_low_high}
\end{figure}

\begin{figure}[h!]
  \centering
  \includegraphics[width=0.7\linewidth]{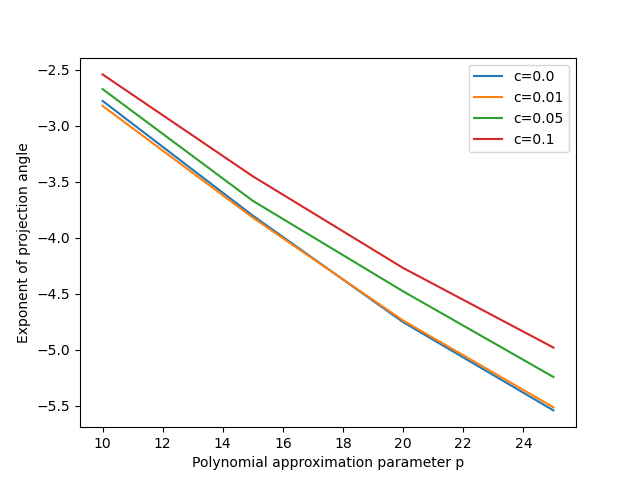}
   \vspace{-6mm}
  \caption{The base-$p$ exponent expressing the error of projection of the force
  corresponding the non-homogeneous problem. The graph features high values of 
  $p$, $10\leq p \leq 25$, and low values of $c$, $0 \leq c \leq 0.1$. This is the
  quantity $s_{nh}(p,c,T)$ in the established notation.}
  \label{fig:exp_force_p_c_non_hom_high_low}
\end{figure}

\begin{figure}[h!]
  \centering
  \includegraphics[width=0.7\linewidth]{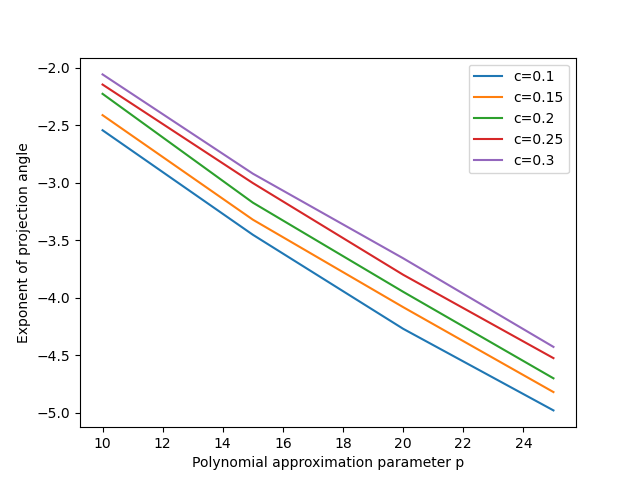}
   \vspace{-6mm}
  \caption{The base-$p$ exponent expressing the error of projection of the force
  corresponding the non-homogeneous problem. The graph features high values of 
  $p$, $10\leq p \leq 25$, and high values of $c$, $0.1 \leq c \leq 0.3$. This is the
  quantity $s_{nh}(p,c,T)$ in the established notation.}
  \label{fig:exp_force_p_c_non_hom_high_high}
\end{figure}

The corresponding quantity for the non-homogeneous problem, $s_{nh} = s_{nh}(p,c,T)$
is also calculated, as defined in eq. \eqref{eqdef sh_nh}, with $c$ switched on.
The exponent $s_{nh} = s_{nh}(p,c,T)$ is plotted in figs.
\ref{fig:exp_force_p_c_non_hom_low_low}, \ref{fig:exp_force_p_c_non_hom_low_high},
\ref{fig:exp_force_p_c_non_hom_high_low}, and \ref{fig:exp_force_p_c_non_hom_high_high},
with the same split for $p$ and for $c$, for the sake of clarity.

The figures show that the effect of very small damping can be in the direction of
increasing the rate of convergence, but the overall tendency is to decrease the rate.
For all $p$ it holds that
\begin{equation}
s_{nh}(p,0.3,T) \leq s_{nh}(p,0,T) + 1.12
\end{equation}
which implies that the approximation of the non-homogeneous problem is less severely
affected by the introduction of damping than the homogeneous one.

\subsubsection{The effect of damping, dependence on $h$}

Finally, the dependence on $h$ is studied, with damping switched on. Once again, the
scaling rule for the projection angle seems impossible to analyze with means other
than a numerical study.

Thus, the exponent $s_{h} = s_{h}(p,c,h)$ is plotted in the following figures, for a
fixed value of $p$ in each one. It is the quantity of eq. \eqref{eqdef sh_h}, which
now reads
\begin{equation} \label{eqdef sh_h_c}
    s_{h} = s_{h}(p,c,h) = \max _{ \max \{ |x_{0}|, |\dx_{0}|\}= 1}
    \frac{\log  \|\pi _{(\dot{\PP} \h)^{\perp}}
    (x_{0}\mathbf{e}_{disp} + \dx_{0}\mathbf{e}_{vel})\|_{L^{2}}
     - s_{h}(p,c,T)\log p}{\log (h/T)}
\end{equation}
with all relevant parameters being variable, so that
\begin{equation}
    \|\pi _{(\dot{\PP} \h)^{\perp}} F\|_{L^{2}} \leq
    ( \max |x_{0}|, |\dx_{0}|) p^{s_{h}(p,c,1)}h^{s_{h}(p,c,h)}, \forall F \in \mathrm{Vec}(\{ \mathbf{e}_{1}
    ,
    \mathbf{e}_{2}\})
\end{equation}

\begin{figure}[h!]
  \centering
  \includegraphics[width=0.7\linewidth]{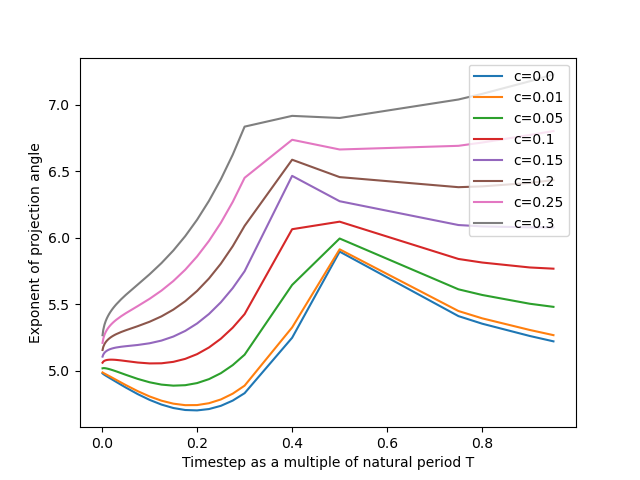}
   \vspace{-6mm}
  \caption{The base-$h$ exponent expressing the error of projection of the force
  corresponding the homogeneous problem for $p=3$. This is the quantity
  $s_{h} = s_{h}(3,c,h)$ in the established notation.}
  \label{fig:exp_force_h_c_hom_3}
\end{figure}

For $p=3$, the exponent is plotted in fig. \ref{fig:exp_force_h_c_hom_3}. Contrary
to what was observed for $h=T$, there is a clear monotonicity, with the exponent
increasing sharply and steadily as $c$ increases. This effect counteracts sufficiently
the loss of convergence due to the deterioration of $s_{h}(p,c,1)$ as $c$ increases,
since, for example,
\begin{equation}
3^{s_{h}(3,0,T)}0.1^{s_{h}(3,0,0.1T)} = 3^{2.255}0.1^{4.780} = \num{1.977e-4}
\end{equation}
while
\begin{equation}
3^{s_{h}(3,0.3,T)}0.1^{s_{h}(3,0.3,0.1T)} = 3^{4.776}0.1^{5.725} = \num{3.579e-4}
\end{equation}
where it is reminded that the value $0.1$ represents the timestep as a fraction of the
natural period of the system. The same calculation for $h=0.01T$ yields, respectively,
$\num{1.312e-4}$ for the undamped system, and $\num{7.677e-4}$ for $c=0.3$.

\begin{figure}[h!]
  \centering
  \includegraphics[width=0.7\linewidth]{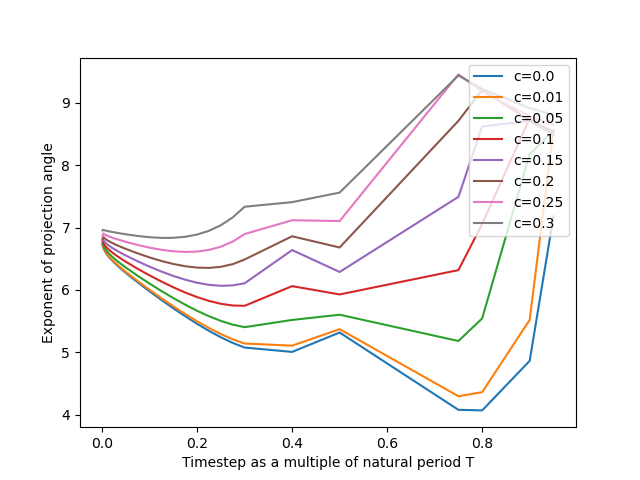}
   \vspace{-6mm}
  \caption{The base-$h$ exponent expressing the error of projection of the force
  corresponding the homogeneous problem for $p=4$. This is the quantity
  $s_{h} = s_{h}(4,c,h)$ in the established notation.}
  \label{fig:exp_force_h_c_hom_4}
\end{figure}

The corresponding graph for $p=4$ can be found in fig. \ref{fig:exp_force_h_c_hom_4},
where the same phenomenon as for $p=3$ can be observed.

\begin{figure}[h!]
  \centering
  \includegraphics[width=0.7\linewidth]{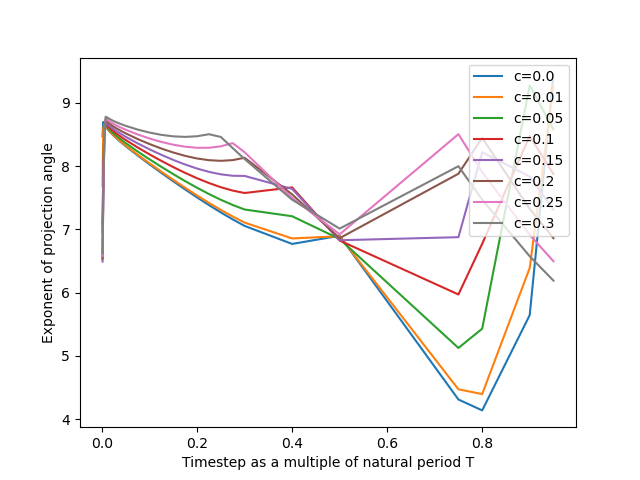}
   \vspace{-6mm}
  \caption{The base-$h$ exponent expressing the error of projection of the force
  corresponding the homogeneous problem for $p=5$. This is the quantity
  $s_{h} = s_{h}(5,c,h)$ in the established notation.}
  \label{fig:exp_force_h_c_hom_5}
\end{figure}

\begin{figure}[h!]
  \centering
  \includegraphics[width=0.7\linewidth]{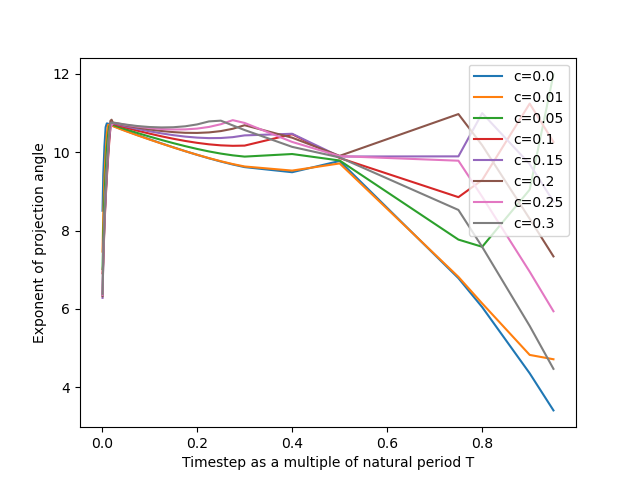}
   \vspace{-6mm}
  \caption{The base-$h$ exponent expressing the error of projection of the force
  corresponding the homogeneous problem for $p=6$. This is the quantity
  $s_{h} = s_{h}(6,c,h)$ in the established notation.}
  \label{fig:exp_force_h_c_hom_6}
\end{figure}

The pattern begins to change for higher order approximations, as can be seen in figs.
\ref{fig:exp_force_h_c_hom_5} and \ref{fig:exp_force_h_c_hom_6}, featuring the curves
for $p=5,6$. The monotonic gains are clear for timesteps up to $\approx 25\%$ of the
natural period of the system, where a tipping point is reached and monotonicity breaks
for stong damping coefficients. The exponent of the timestep remains more favorable,
nonetheless, for damped systems relative to the undamped case.

\begin{figure}[h!]
  \centering
  \includegraphics[width=0.7\linewidth]{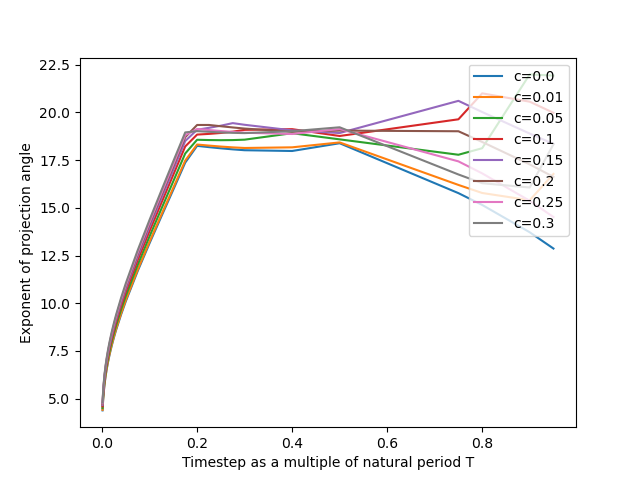}
   \vspace{-6mm}
  \caption{The base-$h$ exponent expressing the error of projection of the force
  corresponding the homogeneous problem for $p=10$ and timesteps smaller than the
  natural period of the system. This is the quantity
  $s_{h} = s_{h}(10,c,h)$ in the established notation.}
  \label{fig:exp_force_h_c_hom_10_low}
\end{figure}

\begin{figure}[h!]
  \centering
  \includegraphics[width=0.7\linewidth]{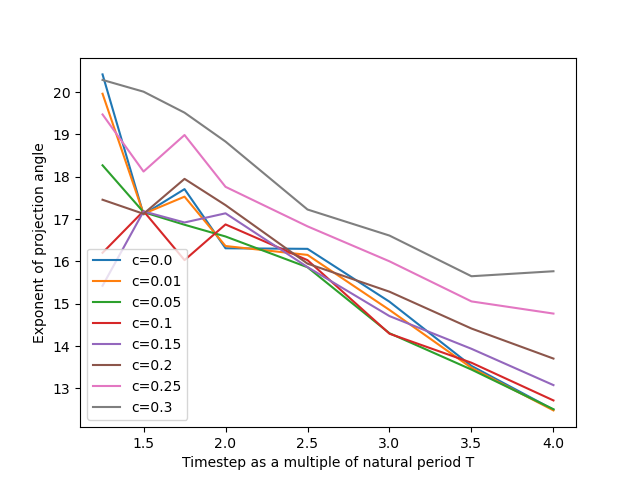}
   \vspace{-6mm}
  \caption{The base-$h$ exponent expressing the error of projection of the force
  corresponding the homogeneous problem for $p=10$ and timesteps greater than the
  natural period of the system. This is the quantity
  $s_{h} = s_{h}(10,c,h)$ in the established notation.}
  \label{fig:exp_force_h_c_hom_10_high}
\end{figure}

Finally, for $p=10$ and timesteps smaller than the natural period the exponent
is plotted in fig.
\ref{fig:exp_force_h_c_hom_10_low}, while for timesteps ranging from $1.5\times$ up to
$4\times T$, in fig. \ref{fig:exp_force_h_c_hom_10_high}.

In the case of timesteps smaller than the natural period of the system, damping leads
in general to gains in convergence speed, while for timesteps greater than the period
the picture is more complicated, as there is no clear monotonicity except for the 
higher part of the interval of damping coefficients studied.

\begin{figure}[h!]
  \centering
  \includegraphics[width=0.7\linewidth]{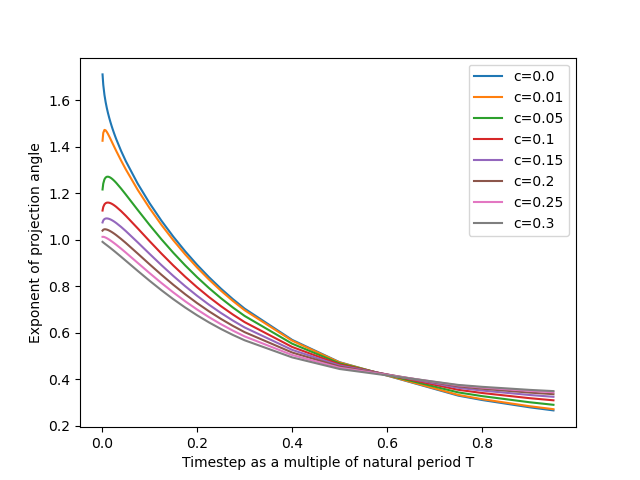}
   \vspace{-6mm}
  \caption{The base-$h$ exponent expressing the error of projection of the force
  corresponding the non-homogeneous problem for $p=3$. This is the quantity
  $s_{nh} = s_{nh}(3,c,h)$ in the established notation.}
  \label{fig:exp_force_h_c_non_hom_3}
\end{figure}

\begin{figure}[h!]
  \centering
  \includegraphics[width=0.7\linewidth]{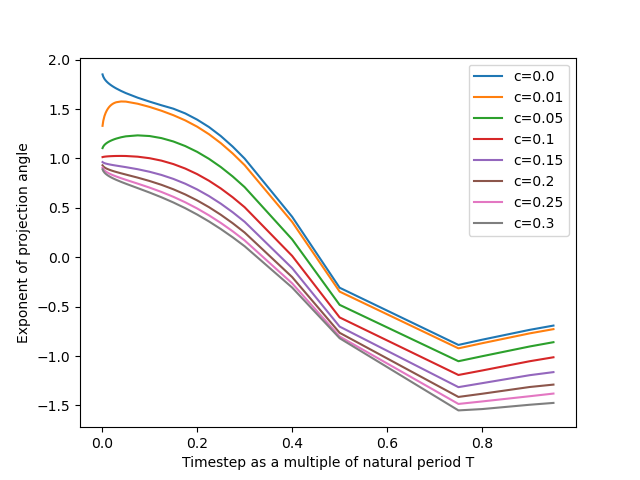}
   \vspace{-6mm}
  \caption{The base-$h$ exponent expressing the error of projection of the force
  corresponding the non-homogeneous problem for $p=4$. This is the quantity
  $s_{nh} = s_{nh}(4,c,h)$ in the established notation.}
  \label{fig:exp_force_h_c_non_hom_4}
\end{figure}

\begin{figure}[h!]
  \centering
  \includegraphics[width=0.7\linewidth]{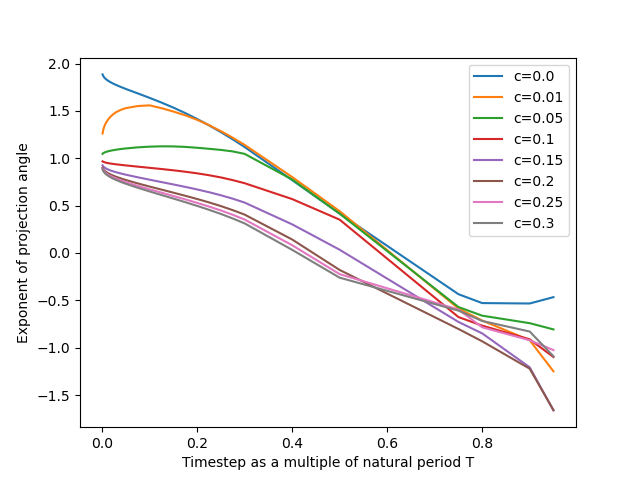}
   \vspace{-6mm}
  \caption{The base-$h$ exponent expressing the error of projection of the force
  corresponding the non-homogeneous problem for $p=5$. This is the quantity
  $s_{nh} = s_{nh}(5,c,h)$ in the established notation.}
  \label{fig:exp_force_h_c_non_hom_5}
\end{figure}

\begin{figure}[h!]
  \centering
  \includegraphics[width=0.7\linewidth]{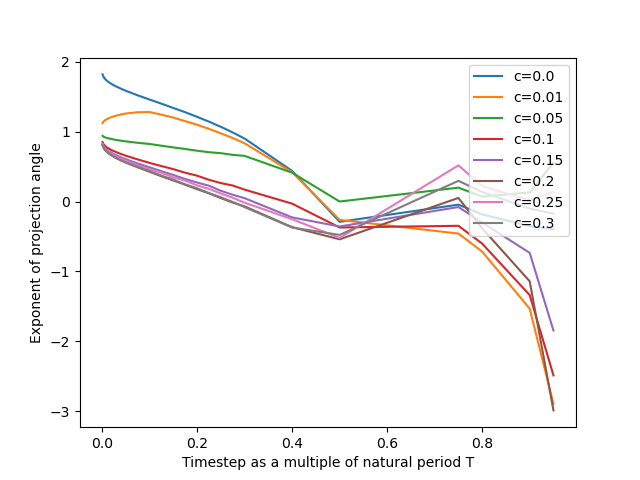}
   \vspace{-6mm}
  \caption{The base-$h$ exponent expressing the error of projection of the force
  corresponding the non-homogeneous problem for $p=6$. This is the quantity
  $s_{nh} = s_{nh}(6,c,h)$ in the established notation.}
  \label{fig:exp_force_h_c_non_hom_6}
\end{figure}

Turning to the non-homogeneous problem now, the quantity
\begin{equation} \label{eqdef sh_nh_c}
    s_{nh} = s_{nh}(p,c,h) = \max _{\substack{ F \in \mathrm{Vec} \{ \mathbf{\tilde{e}}_j \}_{j=3}^{p}
    \\
    \|F \|_{L^{2}} \leq 1}}
    \frac{\log  \|\pi _{(\dot{\PP} \h)^{\perp}} F\|_{L^{2}} - s_{nh}(p,c,T) \log p}{\log h}
\end{equation}
can be calculated, so that
\begin{equation}
    \|\pi _{(\dot{\PP} \h)^{\perp}} F\|_{L^{2}} \leq
    \| F\|_{L^{2}} p^{s_{nh}(p,c,T)} h^{s_{nh}(p,c,h)}, \forall F \in \mathrm{Vec}\{ \mathbf{\tilde{e}}_j \}_{j=3}^{p}
\end{equation}
As before, the results are presented for fixed $p$, varying $c$ and $h$ in figs.
\ref{fig:exp_force_h_c_non_hom_3}, \ref{fig:exp_force_h_c_non_hom_4},
\ref{fig:exp_force_h_c_non_hom_5}, \ref{fig:exp_force_h_c_non_hom_6}.

In this case, the monotonicity is reversed, with convergence speed deteriorating for
all relevant values of the timestep. The dependence of $s_{nh}$ on $c$ becomes milder
as $c$ increases for all relevant values of $h$. This leads to the following,
considerable, deterioration of the convergence rate when passing from $c=0$ to $c=0.3$
for $p=3$ and $h=0.1T$:
\begin{equation}
3^{s_{nh}(3,0,T)}0.1^{s_{nh}(3,0,0.1T)} = 3^{-0.9561223833692502}0.1^{1.1565567840705129} = 0.02439
\end{equation}
while
\begin{equation}
3^{s_{nh}(3,0.3,T)}0.1^{s_{nh}(3,0.3,0.1T)} = 3^{-0.6372893192703653}0.1^{0.8230111827097677} = 0.07463
\end{equation}
For a timestep $h=0.01T$ the quantities read, respectively, $0.0002583$ and
$0.005545$.

\begin{figure}[h!]
  \centering
  \includegraphics[width=0.7\linewidth]{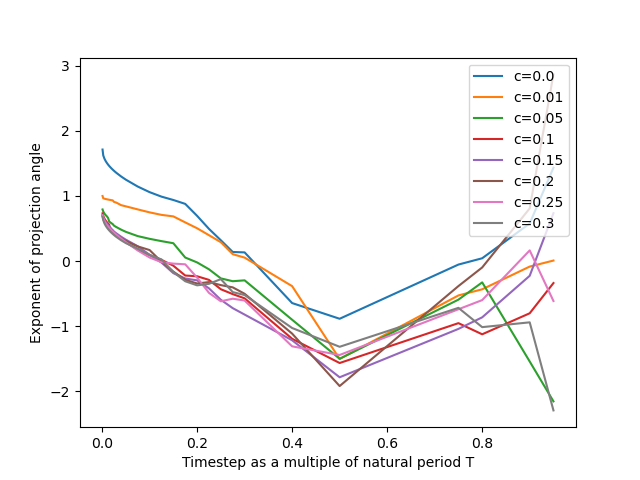}
   \vspace{-6mm}
  \caption{The base-$h$ exponent expressing the error of projection of the force
  corresponding the non-homogeneous problem for $p=10$ and timesteps smaller than the
  natural period of the system. This is the quantity
  $s_{nh} = s_{nh}(10,c,h)$ in the established notation.}
  \label{fig:exp_force_h_c_non_hom_10_low}
\end{figure}

\begin{figure}[h!]
  \centering
  \includegraphics[width=0.7\linewidth]{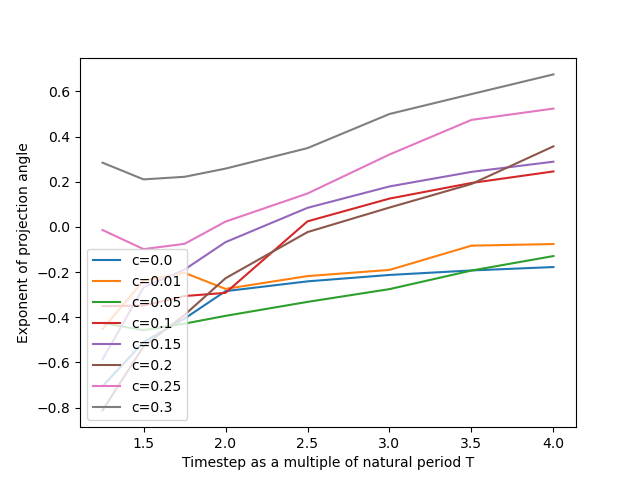}
   \vspace{-6mm}
  \caption{The base-$h$ exponent expressing the error of projection of the force
  corresponding the non-homogeneous problem for $p=10$ and timesteps greater than the
  natural period of the system. This is the quantity
  $s_{nh} = s_{nh}(10,c,h)$ in the established notation.}
  \label{fig:exp_force_h_c_non_hom_10_high}
\end{figure}

Finally, for $p=10$ and timesteps smaller than the natural period of the system,
fig. \ref{fig:exp_force_h_c_non_hom_10_low} shows that in general damping results
in a slower congergence rate, which is practically indifferent to $h$ except for small
values of $h/T$ and $c$. Fig. \ref{fig:exp_force_h_c_non_hom_10_high} shows that
in the case of timesteps greater than the natural period, convergence is also practically
indifferent to $h$, with all exponents close to $0$.

\subsubsection{The aggregate convergence factor} \label{sec conv fac}

In this final section on the error due to the misalignment of spaces, the aggregate
factor is calculated, i.e. the quantity
\begin{equation} \label{eqdef conv factor}
\phi_{\#}(p,c,h) = p^{exp_{\#}(p,c,T)}h^{exp_{\#}(p,c,h)}
\end{equation}
for $\# \in \{ h, nh \} $.
This is necessary due to the convoluted dependence of the exponents on $c$ and $h$,
where it is reminded that the dependence on $p$ is the expected one, i.e. decreasing
in a practically monotonous way. All graphs are logarithmic in order to accommodate
the large range of values of the convergence factor.

\begin{figure}[h!]
  \centering
  \includegraphics[width=0.7\linewidth]{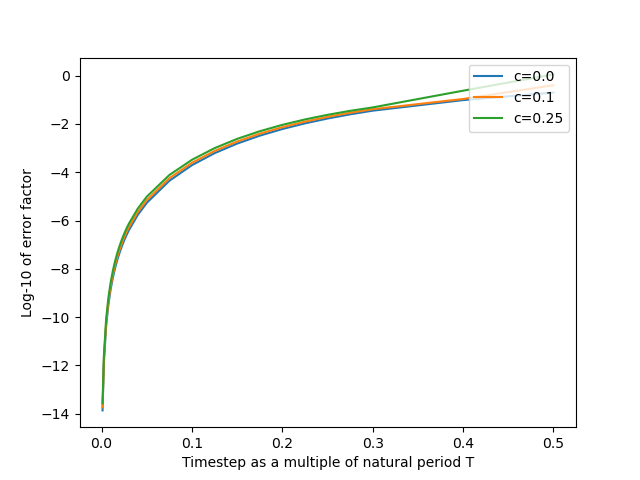}
   \vspace{-6mm}
  \caption{The aggregate factor expressing the error of projection of the force
  corresponding the homogeneous problem for $p=3$ and timesteps smaller than the
  natural period of the system. This is the quantity
  $\phi _{h}(3,c,h)$ in the established notation.}
  \label{fig:force_factor_3_hom_low}
\end{figure}

\begin{figure}[h!]
  \centering
  \includegraphics[width=0.7\linewidth]{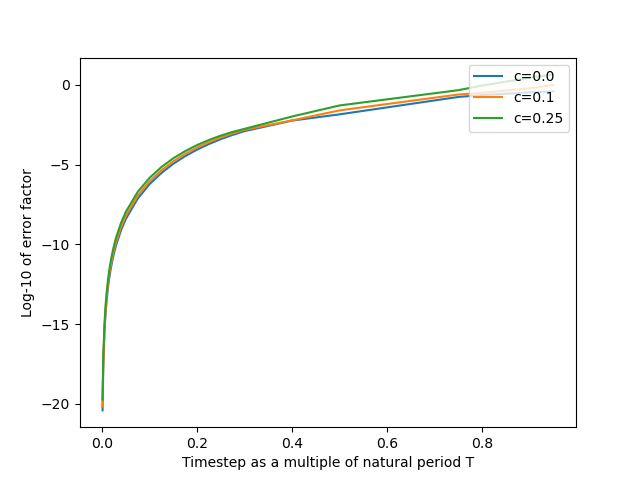}
   \vspace{-6mm}
  \caption{The aggregate factor expressing the error of projection of the force
  corresponding the homogeneous problem for $p=4$ and timesteps smaller than the
  natural period of the system. This is the quantity
  $\phi _{h}(4,c,h)$ in the established notation.}
  \label{fig:force_factor_4_hom_low}
\end{figure}

\begin{figure}[h!]
  \centering
  \includegraphics[width=0.7\linewidth]{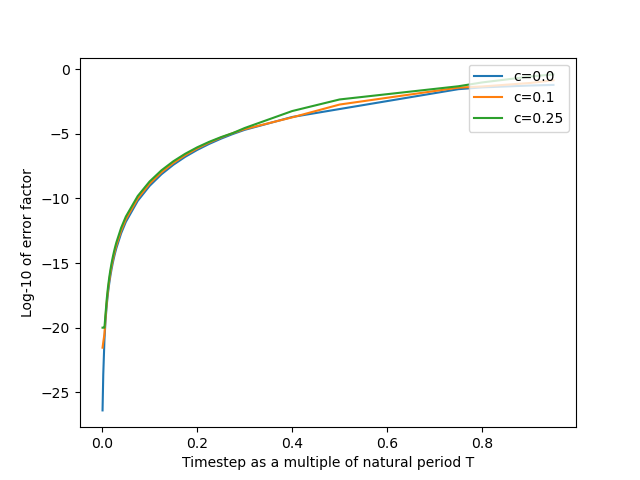}
   \vspace{-6mm}
  \caption{The aggregate factor expressing the error of projection of the force
  corresponding the homogeneous problem for $p=5$ and timesteps smaller than the
  natural period of the system. This is the quantity
  $\phi _{h}(5,c,h)$ in the established notation.}
  \label{fig:force_factor_5_hom_low}
\end{figure}

\begin{figure}[h!]
  \centering
  \includegraphics[width=0.7\linewidth]{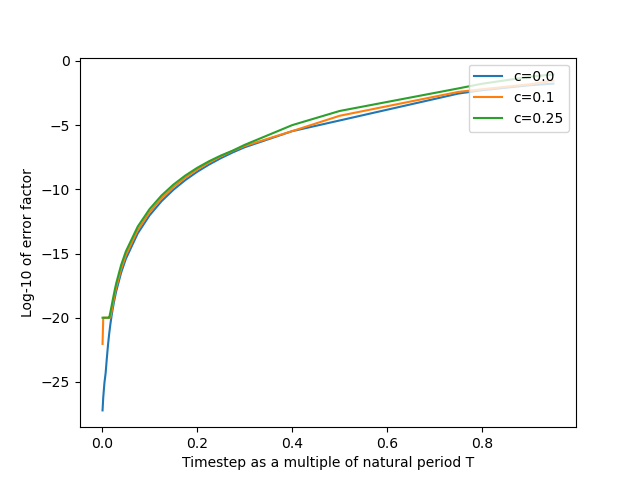}
   \vspace{-6mm}
  \caption{The aggregate factor expressing the error of projection of the force
  corresponding the homogeneous problem for $p=6$ and timesteps smaller than the
  natural period of the system. This is the quantity
  $\phi _{h}(6,c,h)$ in the established notation.}
  \label{fig:force_factor_6_hom_low}
\end{figure}

\begin{figure}[h!]
  \centering
  \includegraphics[width=0.7\linewidth]{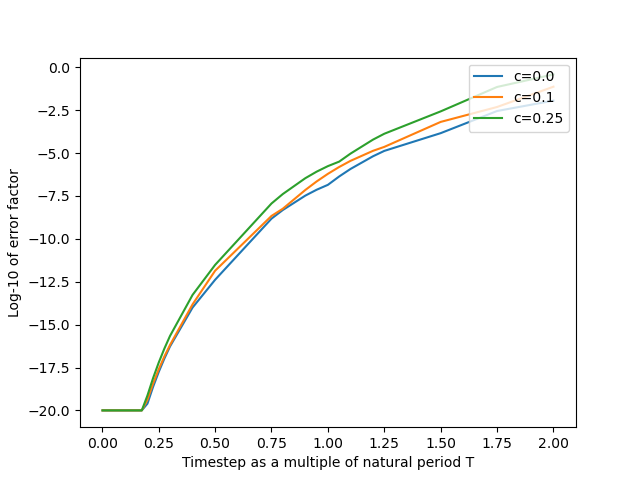}
   \vspace{-6mm}
  \caption{The aggregate factor expressing the error of projection of the force
  corresponding the homogeneous problem for $p=10$ and timesteps smaller than the
  natural period of the system. This is the quantity
  $\phi _{h}(10,c,h)$ in the established notation.}
  \label{fig:force_factor_10_hom_low}
\end{figure}

In the homogeneous case, $\phi_{h}(p,c,h)$ satisfies
\begin{equation}
\| F_{er}\|_{L^{2}} \leq \phi_{h}(p,c,h) \max \{ |x_0|,| \dx_0| \}
\end{equation}
Graphs in figs. \ref{fig:force_factor_3_hom_low}, \ref{fig:force_factor_4_hom_low},
\ref{fig:force_factor_5_hom_low}, \ref{fig:force_factor_6_hom_low} and
\ref{fig:force_factor_10_hom_low}, with a varying range of timesteps. Not all
damping factors have been retained, for clarity, but it can be seen that monotonicity
with respect to $c$ largely holds, with the convergence factor deteriorating as
$c$ increases.

\begin{figure}[h!]
  \centering
  \includegraphics[width=0.7\linewidth]{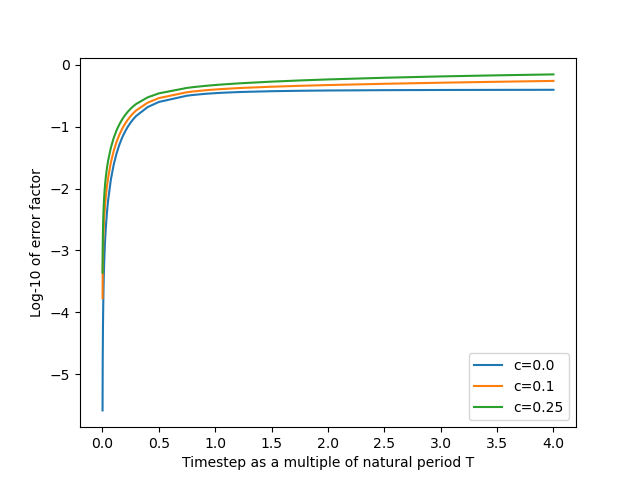}
   \vspace{-6mm}
  \caption{The aggregate factor expressing the error of projection of the force
  corresponding the non-homogeneous problem for $p=3$ and timesteps smaller than the
  natural period of the system. This is the quantity
  $\phi _{nh}(3,c,h)$ in the established notation.}
  \label{fig:force_factor_3_non_hom_low}
\end{figure}

\begin{figure}[h!]
  \centering
  \includegraphics[width=0.7\linewidth]{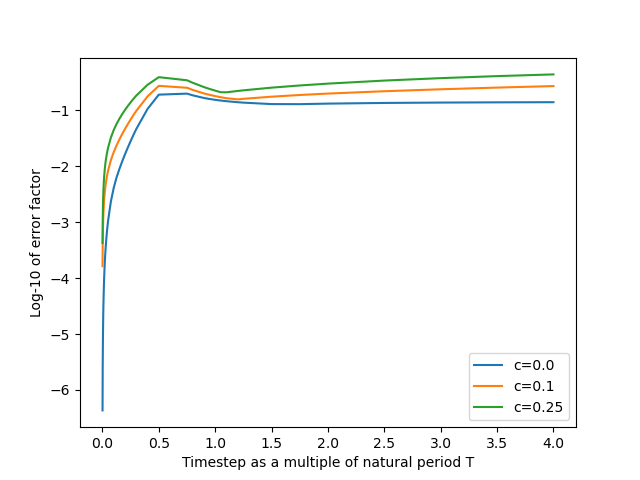}
   \vspace{-6mm}
  \caption{The aggregate factor expressing the error of projection of the force
  corresponding the non-homogeneous problem for $p=4$ and timesteps smaller than the
  natural period of the system. This is the quantity
  $\phi _{nh}(4,c,h)$ in the established notation.}
  \label{fig:force_factor_4_non_hom_low}
\end{figure}

\begin{figure}[h!]
  \centering
  \includegraphics[width=0.7\linewidth]{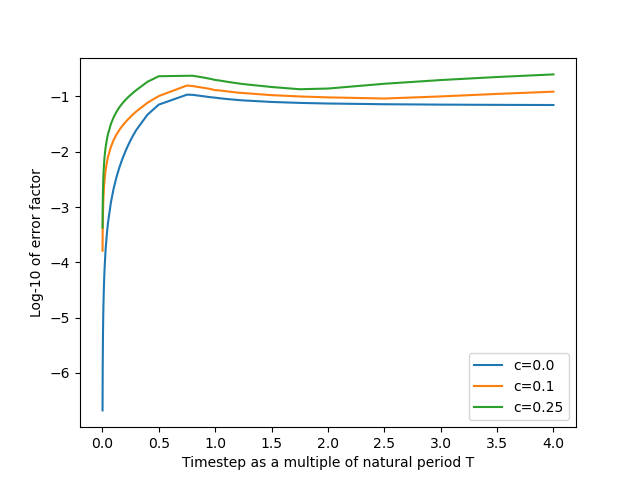}
   \vspace{-6mm}
  \caption{The aggregate factor expressing the error of projection of the force
  corresponding the non-homogeneous problem for $p=5$ and timesteps smaller than the
  natural period of the system. This is the quantity
  $\phi _{nh}(5,c,h)$ in the established notation.}
  \label{fig:force_factor_5_non_hom_low}
\end{figure}

\begin{figure}[h!]
  \centering
  \includegraphics[width=0.7\linewidth]{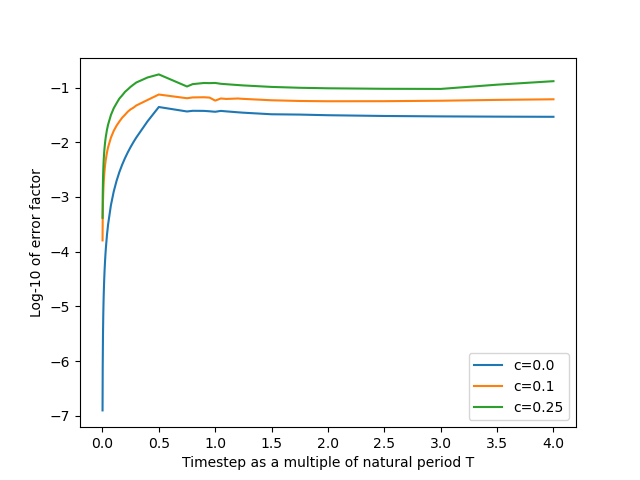}
   \vspace{-6mm}
  \caption{The aggregate factor expressing the error of projection of the force
  corresponding the non-homogeneous problem for $p=6$ and timesteps smaller than the
  natural period of the system. This is the quantity
  $\phi _{nh}(6,c,h)$ in the established notation.}
  \label{fig:force_factor_6_non_hom_low}
\end{figure}

\begin{figure}[h!]
  \centering
  \includegraphics[width=0.7\linewidth]{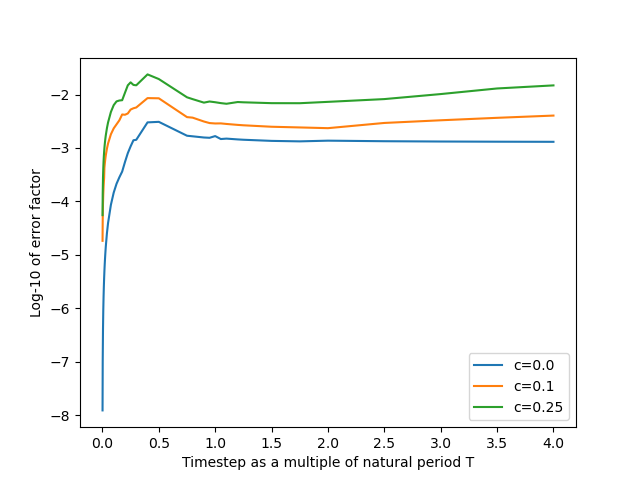}
   \vspace{-6mm}
  \caption{The aggregate factor expressing the error of projection of the force
  corresponding the non-homogeneous problem for $p=10$ and timesteps smaller than the
  natural period of the system. This is the quantity
  $\phi _{nh}(10,c,h)$ in the established notation.}
  \label{fig:force_factor_10_non_hom_low}
\end{figure}

In the non-homogeneous case, $\phi_{nh}(p,c,h)$ satisfies
\begin{equation}
\| F_{er}\|_{L^{2}} \leq \phi_{nh}(p,c,h) \left\Vert \int e^{c\.} \FF (x_{ap})\right\Vert_{L^{2}}
\end{equation}
Graphs in figs. \ref{fig:force_factor_3_non_hom_low},
\ref{fig:force_factor_4_non_hom_low}, \ref{fig:force_factor_5_non_hom_low},
\ref{fig:force_factor_6_non_hom_low} and \ref{fig:force_factor_10_non_hom_low}. The show
monotonicity in $p$ and in $c$, of different directions, and the interesting phenomenon
of indifference of the convergence factor on $h$ for big enough values of the timestep.

\subsection{The error at the end of the timestep} \label{secerror end timestep}

More precise estimates can be be obtained for the error of approximation
at time $t=h$, the timestep of the method, both for the homogeneous and
the non-homogeneous problem, thanks to the estimates obtained in \cite{NKSDOFI},
in prop. 9.9 and 9.12 and the calculations
in the proof of lem. 4.4, and more precisely eq. (38)
and (39). These calculations have as a direct consequence the
estimates of propositions 9.9 and 9.12, i.e.
\begin{equation} \label{eqerror estimate at timestep}
\begin{array}{r@{}l}
| x_{er,h}(h )| &\leq
K
e^{-ch/2} \|F_{er,h}(\.  )\|_{L^{2}_{0}}
\left(
 \|d_{er,h}(\.) \|^{2}_{L^{2}_{0}} +
\|s_{er,h}(\.) \|^{2}_{L^{2}_{0}}
\right)^{1/2}
\\
| \dx_{er,h}(h )| &\leq
e^{-ch /2} \sqrt{h}
\| F_{er,h}(\. )\|_{L^{2}}
+ K |x_{er,h}(h )|
\end{array}
\end{equation}

These expressions are significant because they provide better bounds for the error
propagation from one timestep to the next. This is because the final displacement and
velocity of one timestep are the initial conditions for the next one, and an error in these
initial conditions propagates through the exact solution of the homogeneous SDOF problem.

\begin{prop} \label{propestimates timestep}
    The error of approximation for the displacement and the velocity at
    time $t=h$ satisfy the following estimates
    \begin{equation}
    \begin{array}{r@{}l}
            |x_{er}(h)| &\leq
            K e^{-ch/2} \phi _{h}(p,c,h)
            \left( (p+1)^{-\s} h^{\s} \| F\|_{H^{\s}(0,h)} +
\max \{ |x_{0}|,|\dx_{0}|\} \phi _{h}(p,c,h) + 
 \| F\|_{L^{2}(0,h)} \phi _{nh}(p,c,h) \right)
\\
            |\dx_{er}(h)| &\leq e^{-ch /2} \sqrt{h}
\left( (p+1)^{-\s} h^{\s} \| F\|_{H^{\s}(0,h)} +
\max \{ |x_{0}|,|\dx_{0}|\} \phi _{h}(p,c,h) + 
 \| F\|_{L^{2}(0,h)} \phi _{nh}(p,c,h) \right)
 \\
&\phantom{\max \{ |x_{0}|,|\dx_{0}|\} \phi _{h}(p,c,h) + 
 \| F\|_{L^{2}(0,h)} \phi _{nh}(p,c,h)}
+ K |x_{er,h}(h )|
    \end{array}
    \end{equation}
\end{prop}
In the statement of the proposition, the functions $\phi _{h}$ and $\phi _{nh}$
are defined in eq. \eqref{eqdef conv factor} and have been studied in \S
\ref{sec conv fac}.
\begin{proof}
By prop. 9.7 of \cite{NKSDOFI}, it holds that
$\| x \|_{0} \leq K \| F \|_{0} $ if $F$ is the integral of the force corresponding
to the SDOF problem with homogeneous initial conditions and $x$ is the solution.
Application of this inequality in the case of the force absorbing the initial conditions
and the force corresponding to the error of approximation gives
\begin{equation}
\max \{ \|d_{er,h}(\.) \|_{L^{2}_{0}}, \|s_{er,h}(\.) \|_{L^{2}_{0}} \} \leq
\phi _{h}(p,c,h)
\end{equation}
The remaining factors of the estimate on the final displacement are obtained by direct
application of the definitions of the related quantities on eq.
\eqref{eqerror estimate at timestep}, using cor. 9.6 of \cite{NKSDOFI}.

Regarding the velocity, eq. \eqref{eqerror estimate at timestep} yields directly
the estimate by substitution of the estimate for the displacement and the estimate
for the force used in the estimate of the displacement.
\end{proof}

\subsection{Error propagation}

Proposition \ref{propestimates timestep} concludes the error estimates for the algorithm
proposed herein. This proposition produces the rate of the propagation of
the error, since the final displacement and velocity of timestep $j$ are the initial
conditions of timestep $j+1$, and the error is accumulated by the error of approximation
of the homogeneous problem. To this error, the error due to the projection of the
external force, if the latter is non-zero, is added.

Thus, for initial conditions $x_0$ and $\dx _0$, and external force $f \in H^{-1}$,
the error of displacement at time $t=h$ is given directly by the estimates of prop.
\ref{propestimates timestep}. The errors for times $t \in [0,h]$ is given by
prop. 9.9 of \cite{NKSDOFI}, namely
\begin{equation}
|x_{er}(t)| \leq K  \left( (p+1)^{-\s} h^{\s} \| F\|_{H^{\s}(0,h)} +
\max \{ |x_{0}|,|\dx_{0}|\} \phi _{h}(p,c,h) + 
 \| F\|_{L^{2}(0,h)} \phi _{nh}(p,c,h) \right)
\end{equation}
while for the velocity, by prop. 9.12
\begin{equation}
\|x_{er}\|_{L^{2}(0,h)} \leq K  \left( (p+1)^{-\s} h^{\s} \| F\|_{H^{\s}(0,h)} +
\max \{ |x_{0}|,|\dx_{0}|\} \phi _{h}(p,c,h) + 
 \| F\|_{L^{2}(0,h)} \phi _{nh}(p,c,h) \right)
\end{equation}

If $x_{app,1}$, $\dx_{app,1}$ are the displacement and velocity of the approximate
solution at time $t=h$, then $|x_{er}(2h)|$ is given by the same formula of prop.
\ref{propestimates timestep}, plus a term
\begin{equation}
\max \{ |x_{er}(h)|,|\dx_{er}(h)|\} \phi _{h}(p,c,h)
\end{equation}
due to the error propagation, similarly for the velocity.

Higher regularity for $f$ results in $L^{\infty}$ bounds for the velocity as detailed in
thm. 9.16 of \cite{NKSDOFI}, which can be seen by
taking one derivative of the equation of motion and applying the same estimates on
the velocity instead of the displacement.

\section{Conclusions}

Even though direct comparison with traditional step-wise methods
such as the Newmark method is not easy to make, due to the fundamental difference
in the way the methods are built, the advantage of working with the
weak formulation and then deriving a numerical method can already
be seen in the study of the convergence of the method.

From the estimates it is already clear that, if the order of polynomial
approximation is sufficiently high, timesteps comparable to the
eigenperiod of the system can potentially give competitive results.
This is known to be outside the scope of traditional timestep methods.

Naturally, increasing the order of polynomial approximation comes with
additional computational overhead, but on the other hand the number of
iterations may become smaller without hindering the precision of the
method. The assessment of the trade-off, even though
of crucial importance for applications, goes beyond the scope of the
present work.

The numerical evidence illustrating the quality of the convergence of
the method will be presented in part III of the paper.

\bibliography{aomsample}
\bibliographystyle{aomalpha}

\end{document}